\documentclass[12pt]{amsart}
\usepackage{amsthm,amsmath}
\usepackage{amssymb}
\usepackage{xy}
\usepackage{graphicx}

\newcommand{\de}{\delta}
\newcommand{\ep}{\varepsilon}

\newcommand{\ee}{{\bf e}}
\newcommand\tc[2]{\theta\left[\begin{matrix}#1\\ #2\end{matrix}\right]}

\newcommand{\CC}{{\mathbb{C}}}

\newcommand{\EE}{{\mathbb{E}}}
\newcommand{\PP}{{\mathbb{P}}}

\newcommand{\RR}{{\mathbb{R}}}
\newcommand{\ZZ}{{\mathbb{Z}}}

\newcommand{\calO}{{\mathcal O}}
\newcommand{\calH}{{\mathcal H}}

\newcommand{\calA}{{\mathcal A}}

\newcommand{\calM}{{\mathcal M}}

\newcommand{\calX}{{\mathcal X}}

\newcommand{\op}{\operatorname}
\newcommand{\Sat}{{\calA_g^{\op {Sat}}}}
\newcommand{\Vor}{{\calA_g^{\op {Vor}}}}
\newcommand{\Perf}{{\calA_g^{\op {Perf}}}}

\newcommand{\Sp}{\op{Sp}}

\newcommand{\Pic}{\op{Pic}}

\newcommand{\even}{\op{even}}

\newcommand{\Sing}{\op{Sing}}

\def\diag{\operatorname{diag}}

\theoremstyle{plain}
\newtheorem{thm}{Theorem}[section]
\newtheorem{lm}[thm]{Lemma}
\newtheorem{prop}[thm]{Proposition}
\newtheorem{cor}[thm]{Corollary}

\newtheorem{conj}[thm]{Conjecture}

\theoremstyle{definition}

\newtheorem{rem}[thm]{Remark}

\begin{document}
\title[Semi-abelic varieties of small torus rank]{Principally polarized semi-abelic varieties of small torus rank, and the Andreotti-Mayer loci}
\author{Samuel Grushevsky}
\address{Mathematics Department, Stony Brook University,
Stony Brook, NY 11790-3651, USA}
\email{sam@math.sunysb.edu}
\thanks{Research supported in part by National Science Foundation under the grant DMS-10-53313.}
\author{Klaus Hulek}
\address{Institut f\"ur Algebraische Geometrie, Leibniz Universit\"at Hannover, Welfengarten 1, 30060 Hannover, Germany}
\email{hulek@math.uni-hannover.de}
\thanks{Research is supported in part by DFG grants Hu-337/6-1 and Hu-337/6-2}
\dedicatory{Dedicated to the memory of Eckart Viehweg}

\begin{abstract}
We obtain, by a direct computation,
explicit descriptions of all principally polarized semi-abelic
varieties of torus rank up to $3$. We describe the geometry of their
symmetric theta divisors and obtain explicit formulae for the
involution and its fixed points.


%
These results allow us to give a new proof of the statements about
the dimensions, for small genus, of the loci of ppav with theta
divisor containing two-torsion points of multiplicity three. We also
prove a result about the closure of this set. Our computations are
used in our work \cite{paper} to compute the class of
the closure of the locus of intermediate Jacobians of cubic
threefolds in the Chow and homology ring of the perfect cone compactification of
the moduli space of principally polarized abelian fivefolds.
Our computations will
be used in our forthcoming work \cite{paper} to compute the class of the closure of the locus of intermediate jacobians of cubic threefolds in the Chow ring of the perfect cone compactification of the moduli space of principally polarized abelian fivefolds.
\end{abstract}
\maketitle

\section{Introduction}
The geometry of the moduli space $\calA_g$ of principally polarized abelian varieties (ppav) of dimension $g$ has been studied extensively since the works of Riemann and Jacobi. Still, the structure of the cohomology and Chow rings of $\calA_g$ and its various compactifications remains to a large extent unknown. Much of the recent progress in this direction (and similarly for the moduli of curves $\calM_g$) originated from the idea of studying the tautological subrings of Chow and cohomology. In particular, the tautological ring of $\calA_g$ is the ring generated by the Chern classes of the Hodge vector bundle, and Mumford \cite{mumfordtowards} was perhaps the first to discuss the relations in the tautological ring. Later van der Geer in \cite{vdgeercycles} and Esnault and Viehweg in \cite{esvi} described the tautological ring of $\calA_g$, and of a suitable toroidal compactification (in the cohomology and in the Chow ring, respectively), by proving that the fundamental relation gives all the relations. This fundamental relation is at the heart of many calculations in the Chow ring, and we dedicate this paper to the memory of Eckart Viehweg.

\smallskip
While the tautological ring has been described, it is not clear what the non-tautological classes in Chow and cohomology are. Moreover, for naturally defined geometric subvarieties of $\calA_g$ it is not clear if their classes are tautological. Studying geometrically defined subvarieties of $\calA_g$ could thus shed some further light on the geometry of the space, and this is the subject we pursue here.

Some of the first examples of such subvarieties are the Andreotti-Mayer loci $N_k$. These are the loci of ppav whose theta divisor is singular in dimension at least $k$ :
$$
N_{k,g}:=N_k:= \{(X,\Theta_X) | \dim \Sing (\Theta_X) \geq k\}.
$$
These loci were defined in \cite{anma} as a tool for approaching the Schottky problem --- the question of characterizing Jacobians of curves among all ppav. To this end Andreotti and Mayer proved that the locus of Jacobians is an irreducible component of $N_{g-4}$, and the locus of hyperelliptic Jacobians --- of $N_{g-3}$.

Beauville \cite{beauville} showed that $N_0$ is a divisor, Debarre \cite{debarredecomposes} showed that for $g\ge 4$ the divisor $N_0$ has two irreducible components (one of them the theta-null divisor), and Mumford \cite{mumforddimag} showed that the codimension of $N_1$ is at least 2. In another direction, Ein and Lazarsfeld \cite{eila1} proved the conjecture of Arbarello and De Concini \cite{ardcnovikov} stating that $N_{g-2}$ is equal to the locus of decomposable ppav (those that are products of lower-dimensional ppav).

In general not much is known about the number of components of the Andreotti-Mayer loci, their dimension, or their class in the Chow ring. The study of these loci was recently advanced by Ciliberto and van der Geer, who showed in \cite{civdg1} that
for $1\le k\le g-3$ the codimension of (any component of) $N_k$ is at least $k+2$ (and if $k>g/3$, at least $k+3$). However, this does not seem to be the best bound possible, and they conjecture in \cite{civdg1} that any component of $N_k$ whose general point corresponds to an abelian variety with endomorphism ring $\ZZ$ has codimension at least $\binom{k+2}{2}$. They moreover conjecture that equality only occurs for $g=k+3$ and the hyperelliptic locus, respectively $g=k+4$ and the Jacobi locus. This conjecture remains wide open for $k>1$, and only the case of $k=1$ was fully settled by Ciliberto and van der Geer in \cite{civdg2}. This means that Ciliberto and van der Geer prove that the codimension of $N_1$ is at least 3, this being the first improvement in this direction since the work of Mumford.

The method of their proof is by considering the closure of the locus $N_1$ in a suitable toroidal compactification of $\calA_g$, and studying all degenerations of torus rank up to 2 (we shall discuss this notion in detail below in section \ref{sec:outline}).

It appears that our methods and results on semiabelic varieties of torus rank up to 3 should allow for an extension of their method deeper into the boundary, and in particular would allow one to handle {\it any} subvariety of $\calA_g$ of codimension at most 5, by studying its intersection with the boundary of a toroidal compactification. In particular, we believe our results, together with a suitable interpretation of the results of \cite{civdg2}, should give an approach to the Andreotti-Mayer locus $N_2$. While the conjecture of Ciliberto and van der Geer is that the codimension of $N_2$ is at least 6, our results on semiabelic theta divisors should suffice (once the appropriate generalization of \cite[Prop.~12.1]{civdg2} for other degenerations is obtained --- we will remark on this below) to prove the following weaker statement, which would still be an improvement of the known bound of 4 for the codimension:
%
%
\begin{conj}\label{theo:N2}
The codimension in $\calA_g$ of any component of the locus $N_2$ consisting of simple abelian varieties is at least 5.
\end{conj}
%

\smallskip
Another interesting subvariety of $\calA_g$ is the locus $I^{(g)}$ parameterizing principally polarized abelian varieties
whose theta divisor has a singularity at an odd $2$-torsion point. Casalaina-Martin and Friedman \cite{cmfr,casalaina2} showed that $I^{(5)}$ is equal, within the locus of indecomposable ppav, to the locus of intermediate Jacobians of cubic threefolds. Moreover, Casalaina-Martin and Laza \cite{cmla} studied the closure of $I^{(5)}$, in particular describing the corresponding decomposable ppavs. For arbitrary $g$ the loci $I^{(g)}$ (and generalizations) were investigated in \cite{grsmconjectures}; conjecture 1 in that paper states that $I^{(g)}$, if non-empty, has pure codimension $g$ in $\calA_g$. The
main result of this paper is a new proof by degeneration methods of this conjecture for $g\le 5$:
\begin{thm} \label{theo:Ig}
The locus $I^{(g)}$ is empty for $g\le 2$ and has pure codimension $g$ in $\calA_g$ for $3\le g \leq 5$.
\end{thm}
This result is not new; it follows from the work of Casalaina-Martin and Friedman \cite{cmfr,casalaina2} characterizing $I^{(5)}$. However, their results were obtained using the theory of Prym varieties, and our methods are completely different.

Our
%
next
result is what we need for computing the class of the locus of intermediate Jacobians in \cite{paper}. To explain its statement, we need to discuss the compactifications of $\calA_g$. Indeed, the main idea of the proof of the theorems above is to study and understand the closures $\overline{N_2}$ and $\overline{I^{(g)}}$ of the loci $N_2$ and $I^{(g)}$ in a suitable (toroidal)
compactification of $\calA_g$. There are two toroidal compactifications which are of interest to us: the perfect
cone compactification $\Perf$ and the second Voronoi cone compactification $\Vor$. Both are toroidal compactifications, i.e.~are given by a suitable fan (or rational polyhedral decomposition)
in the rational closure of the cone of real semi-positive $g \times g$ matrices, namely the
perfect cone and the second Voronoi decomposition respectively. Both $\Perf$ and $\Vor$ have an intrinsic
geometric meaning: the first is the minimal model of $\calA_g$ for $g \geq 12$, in the sense of the minimal model program (this
was proven by Shepherd-Barron \cite{shepherdbarron}), the latter represents (at least up to normalization and up to extra components) a geometrically meaningful functor, namely that of stable polarized semi-abelic varieties (this was proven by Alexeev \cite{alexeev}). For $g\leq 5$, but not in general, $\Vor$ is a blow-up of $\Perf$.

Analytically, the locus $I^{(g)}$ can be described as the zero locus of the gradients $f_m$ of the theta function at odd two-torsion points $m$ (see next section for a precise definition). It turns out that $f_m$ extend to sections $\tilde {f_m}$ of some vector bundle on (a level cover of) $\Perf$. We denote by $\overline{G^{(g)}}$ the (projection to $\Perf$ from a level cover of the) zero locus of $\tilde{f}_m$.
We then have the obvious inclusion $\overline{I^{(g)}} \subset \overline{G^{(g)}}$,
as $\overline{G^{(g)}}$ is the zero locus of some holomorphic functions
vanishing identically on $I^{(g)}$ (and thus also on $\overline{I^{(g)}}$). We note that since locally $I^{(g)}$ is given by the vanishing of each of the $g$ components of the gradient vector $f_m$, the codimension of any component of $I^{(g)}$ in $\calA_g$ is at most equal to $g$.
%
The following
result will be crucial in our forthcoming paper \cite{paper}, where it will be used to compute the classes of the loci $\overline{I^{(g)}}$ in the Chow rings of $\Perf$ for $g\le 5$.
\begin{thm}\label{theo:Gg}
The locus $\overline{G^{(g)}}\subset\Perf$ has no irreducible components contained in $\partial\Perf$ of codimension less than 6 in $\Perf$. In particular, for $g\le 5$, we have $\overline{I^{(g)}}=\overline{G^{(g)}}$.
\end{thm}
To see that the first statement of the theorem implies the second, note that for $g\le 5$ the locus $\overline{G^{(g)}}$, being of codimension at most $5$, has no irreducible components contained in the boundary, and thus any of its irreducible components is a closure of a component of $I^{(g)}$.

\medskip
Besides these main theorems, the core of the paper (and of the proofs for all of the above) is an explicit description of the geometry of various types of semi-abelic varieties --- which is really our main result. We are essentially able to extend the degeneration computations originated by Mumford for the partial compactification two steps further into the boundary.

In order to study the closures of the subvarieties $N_k$ and $I^{(g)}$ in $\Perf$ or $\Vor$ it is
necessary to understand the geometry of the semi-abelic varieties and their (semi-abelic polarization) theta divisors, associated to boundary
points of $\Vor$. Recall that every toroidal compactification of $\calA_g$ maps to the Satake compactification
$$
\Sat= \calA_g \sqcup \calA_{g-1} \sqcup \calA_{g-2} \ldots \sqcup \calA_0.
$$
We denote the projections  $P: \Perf \to \Sat$ and $Q: \Vor \to\Sat$.
We denote $\beta_k^0:= P^{-1}(\calA_{g-k})$ the open strata of the perfect cone compactification, and denote by
$\beta_k:=\beta_k^0\sqcup\ldots\sqcup\beta_g^0=P^{-1}(
\calA_{g-k}^{\rm Sat})$ the closed strata. We say that a semi-abelic variety has torus rank $k$ if the corresponding point of $\Perf$ lies in $\beta_k^0$ (resp.~for $\Vor$, lies in $Q^{-1}(\calA_{g-k}^{\op{Sat}}$)

Since the second Voronoi and the perfect cone decomposition coincide up to
genus $3$ it follows in particular that
$\Perf \setminus P^{-1}(\calA_{g-4}^{\operatorname {Sat}}) \cong \Vor \setminus Q^{-1}(\calA_{g-4}^{\operatorname {Sat}})$.
In fact more is true as we shall discuss in Section \ref{sec:outline}: the two compactifications coincide outside codimension
at least $6$ in $\Perf$.
We shall make use of this fact.
One of the results of this paper is a complete description of the geometry of semi-abelic varieties of torus rank at most $3$ and their theta divisors, describing explicitly their geometry, including the geometry of the
polarization divisor. One of the main applications of this description is to computing the class of the closure of the locus $I^{(g)}$ in the Chow ring of $\Perf$, which we will accomplish in a forthcoming paper \cite{paper}. For this reason we also describe explicitly
the involutions on the semi-abelic varieties, their fixed points, and the gradients of the (semi-abelic) theta function
at these fixed points.
We shall derive these descriptions from toric geometry together with some simple properties of
the (symmetric) theta divisor which follow from general principles. Alternatively one can also make Alexeev's construction
explicit. We chose the geometric approach as the more elementary method, as this gives us the geometric
properties which we require.

\section*{Acknowledgements}
We would like to thank Gerard van der Geer for describing the problem to us, and for discussing his work on the subject.
We are grateful to Valery Alexeev for discussions and explanations regarding degenerations of abelian varieties
and their polarizations.
We are greatly indebted to Frank Vallentin, who has patiently answered our numerous questions on cone decompositions.
Finally we thank Niko Behrens for producing the pictures used in this paper. The first author thanks Leibniz Universit\"at Hannover for hospitality during his visits there when some of this work was done.

\section{Theta functions and their gradients}
We denote by $\calH_g$ the Siegel upper half plane of genus $g$ and consider the quotient $\calA_g = \calH_g/\Sp(2g,\ZZ)$  --- the
moduli space of principally polarized abelian varieties (ppav). The Riemann theta function
$\theta(\tau,z)$ is a function of $\tau\in\calH_g$ and $z\in\CC^g$ and is given by
$$
 \theta(\tau,z):=\sum\limits_{n\in\ZZ^g}\ee (\pi i n^t\tau n +2\pi i n^tz),
$$
where for future use we denote $\ee(x):=\exp(2\pi i x)$ the exponential function.
The theta function is quasiperiodic and satisfies
\begin{equation}\label{qperiodic}
 \theta(\tau,z+n\tau+m)=\ee(-n^t\tau n/2-n^tz)\theta(\tau,z)
\end{equation}
for any $n,m\in\ZZ^g$.

For an abelian variety $A$, we denote by $A[2]$ the set of two-torsion points on it; as a group,
$A[2]\cong (\ZZ/2\ZZ)^{2g}$.
Analytically the points in $A[2]$ can be labeled  $m=(\tau\ep+\de)/2$, for $\tau$ being the period matrix of $A$, and
$\ep,\de\in (\ZZ/2\ZZ)^g$. We denote $\sigma(m):=\ep\cdot\de\in\ZZ/2\ZZ$ and call it the {\it parity} of $m$.
Accordingly we
call $m$ even or odd depending on whether $\sigma(m)$ is $0$ or $1$, respectively.
This is equivalent to the point $m$ not lying (resp.~lying) on the theta divisor for a
generic $\tau$ (i.e.~for a two-torsion point $m$ the function $\theta(\tau,m)$ is identically zero if and only if $m$ is odd).

For a point $m=(\tau\ep+\de)/2$ we denote the theta function with (half-integer) characteristic
$$
 \tc\ep\de(\tau,z):=\theta_m(\tau,z):=
$$
$$
=\sum\limits_{n\in\ZZ^g}\ee ((n+\ep/2)^t\tau (n+\ep/2)/2
 + (n+\ep/2)^t(z+\de/2)).
$$
As a function of $z$, the theta function with characteristic is even or odd depending on the parity of the characteristic.
In particular, all odd theta constants with characteristics (the values of theta functions with characteristics at $z=0$)
vanish identically.

We recall the level subgroups of $\Gamma_g:=\Sp(2g,\ZZ)$ defined as follows:
$$
 \Gamma_g(n):=\left\lbrace \gamma=\begin{pmatrix} A&B\\ C&D\end{pmatrix}\in\Gamma_g \right|\left. \gamma\equiv \begin{pmatrix} 1&0\\ 0&1\end{pmatrix}\mod n\right\rbrace
$$
$$
 \Gamma_g(n,2n):=\left\lbrace\gamma\in\Gamma_g(n)\mid \diag(A^tB)\equiv\diag(C^tD)\equiv 0
 \mod 2n\right\rbrace.
$$
The moduli space of ppav is then $\calA_g=\calH_g/\Gamma_g$, while the level moduli spaces $\calA_g(n):=\calH_g/\Gamma_g(n)$ and $\calA_g(n,2n):=\calH_g/\Gamma_g(n,2n)$ are finite covers of $\calA_g$.

We can thus compute for an odd point $m\in A[2]$
$$\begin{aligned}
 f_m(\tau)&:={\rm grad}_z\theta(\tau,z)|_{z=m}={\rm grad}_z\theta(\tau,z+(\tau\ep+\de)/2)_{z=0}\\
&\ = \ee(-\ep^t\tau\ep/8-\ep^t\de/4-\ep^tz/2){\rm grad}_z\tc\ep\de(\tau,z)|_{z=0}
\end{aligned}
$$
since $\tc\ep\de(\tau,0)=0$ for an odd two-torsion point $m$. It is known that this gradient is a vector-valued modular form for $\Gamma_g(4,8)$ for the
representation $\det^{\otimes 1/2}\otimes{\rm std}:\Gamma_g\to GL(g,\CC)$. In other words, we have
\begin{equation}\label{fmbundle}
 f_m\in H^0(\calA_g(4,8),\det\EE^{\otimes 1/2}\otimes\EE),
\end{equation}
where $\EE$ is the Hodge vector bundle on $\calA_g$. Note that this bundle extends to any toroidal compactification of $\calA_g$, see \cite{mumhirz}.

The group $\Gamma_g(2)/\Gamma_g(4,8)$ acts on the gradients by certain signs, and thus the zero locus
\begin{equation}\label{Gmdef}
 G_m:=G_{\epsilon,\delta}
 =\{\tau| f_m(\tau)=0\}
\end{equation}
is a well-defined subvariety of $\calA_g(2)$. We refer to \cite{grsmodd1,grsmconjectures} for a
more detailed discussion of the properties of the
gradients of the theta function and further questions on loci of ppavs with points of high multiplicity on the theta divisor.

Finally we denote by
$$
 I^{(g)}:=p(G_m)\subset\calA_g
$$
the locus of ppavs for which {\it some} $f_m$ vanishes --- note that it follows from the fact that $\Gamma_g$ permutes the $f_m$ that the
projection of $G_m$ to $\calA_g$ does not depend on $m$. For further use we also denote by $\overline{I^{(g)}}\subset\Perf$ the closure in the perfect cone compactification.
The sections $f_m$ extend to sections of a suitable vector bundle on the perfect cone decomposition $\Perf$.
This is most easily described on the level $8$ cover $\Perf(8)$. We shall denote the boundary components of this
compactification by $D_n$. Recall that these are enumerated by the index $\pm n \in (\ZZ / 8\ZZ)^{2g}.$ We denote
the mod $2$ reduction of $n$ by $n_2$. Moreover we call an element $k=(k_1,k_2) \in (\ZZ / 2\ZZ)^{2g}$
{\em even} or {\em odd} if the product of the length $g$ vectors $k_1 \cdot k_2$ is even or odd respectively.
In \cite{paper} we show the following result:
\begin{prop} \label{prop:extension}
The sections $f_m$ extend to global sections
$$
 \tilde f_m\in H^0\left(\Perf(8),\det(\EE)^{\otimes(1/2)}\otimes\EE\otimes \calO\left(
-\sum\limits_{\lbrace n| m+n_2\rm \even\rbrace/n\sim-n}D_n\right)\right),
$$
which do not vanish at a general point of any irreducible component of the boundary divisor.
\end{prop}
\begin{proof}
This follows from a computation of the Fourier expansion of the sections $f_m$. For details see \cite[Section 3]{paper} and also the computations in \cite{grsmconjectures}.
\end{proof}
\begin{rem}
Note that the statement of the above proposition includes the claim
that $\sum\limits_{\lbrace n| m+n_2\rm \even\rbrace/n\sim-n}D_n$ is a Cartier divisor, not just a
Weil divisor, which is not immediate as
$\Perf(8)$ is not a smooth variety for $g \geq 4$. However this was shown in \cite[Section 3]{paper}, and also follows from the computations in \cite{grsmconjectures}.
\end{rem}
We denote by $\overline {G_m}\subset\Perf(4,8)$ the vanishing locus of $\tilde{f}_m$. We note that by this definition we automatically have $\overline{ I^{(g)}}\subset p(\overline{G_m})$, and thus if we show that $\tilde{f}_m$ does not vanish in codimension less than $g$ for $g\le 5$, this would suffice to prove theorem \ref{theo:Ig}.

\section{Outline of approach} \label{sec:outline}
Throughout the paper we shall make extensive use of the description of $\Vor$ as (the normalization of) a moduli space
of polarized semi-abelic varieties. The reader is referred to \cite{alna},\cite{alexeev}, \cite{olsson} where details can be
found. Here we shall recall the facts which we will require for our purposes and outline the approach which we shall
take.

Recall that Alexeev \cite{alexeev} showed that (possibly up to normalization) the
Voronoi compactification $\Vor$ represents a functor. The geometric
objects are pairs $(X,\Theta)$ where $X$ is a stable semi-abelic variety (in particular it admits an action of a
semi-abelian variety) and $\Theta \subset X$ is a theta divisor (for a precise definition
see \cite[Def.~1.1.9]{alexeev}). For smooth objects
the functors of principally polarized abelian varieties and that of stable polarized semi-abelic varieties are
isomorphic: an abelian variety is a torsor over itself and $\Theta$ is just the classical theta divisor (which is unique
up to translation).

The boundary strata of any toroidal compactification correspond to orbits of cones in the corresponding rational
polyhedral cone decomposition (fan). To every such
orbit corresponds a type of polarized semi-abelic varieties. If the relative interior of a given cone contains
symmetric matrices of rank $k$, then the corresponding semi-abelic varieties have torus rank $k$, i.e.~the normalization
of the semi-abelic variety
is fibered over an abelian variety $B$ of dimension $g-k$. The semi-abelic variety will in
general have several components, each fibering over $B$ with fibers being toric varieties of dimension $k$.
The precise nature of these components, as well as the combinatorics of
their gluings, is determined by a periodic decomposition of $\RR^k$ into polytopes, determined by the cone in question,
which is invariant under the action of
$\ZZ^k$.

In \cite{alna},\cite{alexeev},\cite{olsson} the general theory of degeneration
data is described, from which it is also possible to obtain the description of the theta divisor.
In particular, Valery Alexeev explained to us how to easily derive the formulae for the theta divisor of standard semi-abelic varieties (see below for a definition) directly from this general framework.
Note, however, that in this general approach the calculations dealing with reducible semi-abelic
varieties with non-trivial abelian part are already quite complicated, while our approach
is direct and sufficient for our purposes here,
and also provides an explicit and direct confirmation of the general framework. We also compute the involution and study its fixed points explicitly.

\smallskip
In this paper we are mostly interested in the perfect cone decomposition $\Perf$.
The perfect cone decomposition and
the second Voronoi decomposition coincide in genus $g \leq 3$. For $g=4,5$ the second Voronoi decomposition is a refinement
of the perfect cone decomposition \cite{ryba}, \cite{erry2}, i.e.~$\Vor$ is a blow-up of $\Perf$ for these genera.
For $g=4$ the situation is very simple: $\calA^{\operatorname {Vor}}_4$ is the blow-up of $\calA^{\operatorname {Perf}}_4$
in one point. This corresponds to a subdivision of a non-basic $10$-dimensional cone, the so-called second perfect
cone. In $g=5$ the situation is more complicated.
For $g \leq 5$ we let $Z_g \subset \Perf$ be the center of the modification $ \sigma: \Vor \to \Perf$. Then $Z_g = \emptyset$ for
$g \leq 3$ and $Z_4$ is a point.
We shall use the fact that the codimension of $Z_g$ in $\Perf$ is greater than $5$. To see this it is enough to
consider the part of $\Perf$ which lies over $\calA_g \sqcup \calA_{g-1} \ldots \sqcup \calA_{g-5}$ as the complement
of this set in $\Perf$ has codimension $6$. In other words we must consider all those cones in the perfect cone decomposition
which are not also cones in the second Voronoi decomposition. There are no such cones for $g \leq 3$. We already mentioned
that there is precisely one such cone for $g=4$ and this is in codimension $10$. Finally there is only one
stratum in $\Perf$ over $\calA_{g-5}$ of codimension $5$ and this is the stratum belonging to the cone
$\langle x_1^2, \ldots x_5^1 \rangle $. But this cone is a face of the so-called principal cone and is thus
also contained in the second Voronoi decomposition.

A partial compactification $\Perf \setminus Z_g \cong \Vor \setminus \sigma^{-1}(Z_g)$
can be constructed by adding all those strata to $\calA_g$ which correspond to cones belonging to both
the perfect cone and the second Voronoi decomposition. Below we list all cones which correspond to strata
of codimension $\leq 5$ in $\Perf$.
All of these cones belong to both the perfect cone and the second Voronoi decomposition.
These are the strata which we will have to handle.

To explain the geometric situation which we will encounter we start with the simplest case of semi-abelic varieties, the
so-called {\it standard degenerations}. We shall treat these in detail
in the next section.
The starting point for such a standard degeneration is the decomposition of $\RR^k$ into unit cubes.  In this case the normalization of the
semi-abelic varieties is
irreducible, and is a fibration over some $B\in\calA_{g-k}$ with fibers isomorphic
to $(\PP^1)^k$, the toric variety corresponding to a $k$-dimensional cube (we refer to \cite{fultontoric} for a
description of the toric variety associated to a polytope).
All other decompositions which we shall have to consider are
subdivisions of these cubes.

The action of $\ZZ^k$ by translations on the periodic decomposition of $\RR^k$ into polytopes also tells us
how to {\it glue} the normalization to get the semi-abelic variety
(i.e.~what identifications have to be made to obtain the semi-abelic variety from its normalization).
In the case of the
cube these are simply identifications of each pair of the opposite faces. These identifications are {\it gluing
maps} of toric (sub)bundles over $B$. The gluing maps fibers to fibers, and thus induces a gluing map on
the base of the fibration. Such a gluing map $B\to B$ then is a translation by some point in $B$, which
we then call the {\it shift}. Once the shift is determined, we need to understand the map of the toric variety
fiber to another such fiber. Such a map induces a map of the corresponding polytope, and there is
the discrete data --- the combinatorics of which faces of the polytope are glued to which --- and the
continuous parameters which are the multiplicative parameters for the map $(\CC^*)^k\to(\CC^*)^k$. The shift and
continuous gluing parameters account for the moduli of semi-abelic varieties of a given combinatorial type (corresponding to a give cone, i.e.~for given toric variety fiber).
The next step we take is to
normalize (simplify) these continuous parameters by suitable choices of coordinates on each toric variety, wherever possible.

The polytope picture also gives us a description of the theta divisor on the toric part. On each irreducible component of the normalization of the semi-abelic variety we know the numerical class of the
polarization divisor, and we will thus start by writing down a general divisor in that numerical class.
Since the theta divisor on the normalization
is the lifting of the divisor on the semi-abelic variety, its restrictions to the glued faces must be glued, which
gives some restrictions on the formula for the semi-abelic theta divisor.
Finally, the involution on the toric part corresponds to reflection with respect
to the center of the cube: notice that since the cube is dissected, the involution may permute the irreducible
components of the normalization of the semi-abelic variety. Requiring the theta divisor to be preserved
by the involution on the semi-abelic variety will allow us to fix the remaining free shift and gluing parameters.

In the following table we list all the strata in codimension at most $5$ in $\Perf$ (notice that the list
for $\Vor$ would be considerably longer, we refer the reader to \cite{erry2} and \cite[Ch.~4]{vallentinthesis}
for more details on the structure of the
perfect cone and second Voronoi compactifications). We list generators of
a cone in each orbit, describe the slicing of the standard cube defined by this cone and explain
which toric varieties are associated to these polytopes. In this table $F(2,2)$ stands for the toric variety
associated to an octahedron: it is the intersection of $2$ quadrics $\{y_0y_1=y_2y_3=y_4y_5\}$
in $\PP^5$ (the corresponding cone is generated by $x_1^2,x_2^2,x_3^2,(x_1+x_2+x_3)^2$ or
by the ones given in the table) . We also denote $F_4$ the cone over a smooth quadric surface, given by
$\{y_1y_2=y_3y_4\}$ in $\PP^4$ --- it is the toric variety whose associated polytope is the square pyramid.
Finally we denote by $X$ the toric variety whose associated polytope is a 4-dimensional cube with two 4-dimensional
simplices cut off at the opposite corners (this cone is
equivalent to the one generated by $x_1^2,x_2^2,x_3^2,x_4^2,(x_1+x_2+x_3+x_4)^2$,
but the one below is more useful for our computations) --- its geometry is unimportant to us.

\begin{table}\label{tab:listcones}\small
\begin{tabular}{|c|c|c|c|c|}
\hline
&&&&\\
${\Sat}$& Forms generating  & Codim & Type of & Toric part of the\\
& the cone & in $\Perf$ & polytope & semi-abelic variety\\
&&&&\\
\hline
&&&&\\
$\beta_1$& $x_1^2$ &1& interval & $\PP^1$\\
&&&&\\
\hline
&&&&\\
$\beta_2$& $x_1^2,x_2^2$&2& square & $\PP^1\times\PP^1$\\
$\beta_2$& $x_1^2,x_2^2,(x_1-x_2)^2$&3& 2 triangles & $\PP^2\sqcup\PP^2$ \\
&&&&\\
\hline
&&&&\\
$\beta_3$& $x_1^2,x_2^2,x_3^2$&3& cube  & $\PP^1\times\PP^1\times\PP^1$\\
$\beta_3$& $x_1^2,x_2^2,(x_1-x_2)^2,x_3^2$ &4&  2 prisms & $\PP^1\times\PP^2\sqcup\PP^1\times\PP^2$ \\
$\beta_3$&$x_1^2,x_2^2,$&4&1 octahedron& $F(2,2)\sqcup 2\PP^3$\\
 &$(x_1-x_3)^2,(x_2-x_3)^2$ & & and 2 tetrahedra&\\
$\beta_3$&$x_1^2,x_2^2,x_3^2,$&5&2 pyramids&  $2F_4\sqcup2\PP^3$\\
 & $(x_1-x_3)^2,(x_2-x_3)^2$& & and 2 tetrahedra&\\
$\beta_3$&$x_1^2,x_2^2,x_3^2,(x_1-x_2)^2$&6&6 tetrahedra&$6\PP^3$\\
&$(x_1-x_3)^2,(x_2-x_3)^2$&&&\\
\hline
&&&&\\
$\beta_4$&$x_1^2,x_2^2,x_3^2,x_4^2$&4&$4$-dim cube& $(\PP^1)^{\times 4}$\\
$\beta_4$&$x_1^2,x_2^2,(x_1-x_2)^2,$&5&product of& $\sqcup_2(\PP^1)^{\times 2}\times\PP^2$\\
 & $x_3^2,x_4^2$& &square and &\\
 & & & 2 triangles&\\
$\beta_4$&$x_1^2,x_2^2,x_3^2,$&5& product of&$\PP^1\times(F(2,2)\sqcup 2\PP^3)$\\
 &$(x_1-x_3)^2,(x_2-x_3)^2,$ & & interval with&\\
 & $x_4^2$& & 1 octahedron&\\
  & & & and 2 tetrahedra&\\
$\beta_4$&$x_1^2, x_2^2, (x_1-x_4)^2$&5&two 4-dim simplices  &$X\sqcup 2\PP^4$\\
&$(x_2-x_3)^2,(x_3-x_4)^2$&&and a 4-dim cube & \\
&&&with simplices cut off&\\
$\beta_4$&\dots&10&24 4-dim simplices&$24\PP^4$\\
\hline
&&&&\\
$\beta_5$&$x_1^2,x_2^2,x_3^2,x_4^2,x_5^2$&5&$5$-dim cube& $(\PP^1)^{\times 5}$\\
&&&&\\
\hline
\end{tabular}
\end{table}

This table gives explicit geometric descriptions of all combinatorial types of principally polarized semi-abelic varieties of torus rank up to 3, and also of all those that give loci of codimension at most 5 in $\Perf$. We would also like to point out that Nakamura \cite{nakamura} has given an explicit description of all semi-abelic varieties
in genus $3$, in a slightly different language.

For each of the cases in this table we will compute explicitly the polarization (semi-abelic theta) divisor, the involution, and the gradients of the semi-abelic theta divisor at the fixed points of the involution. This analysis is also essential for our paper \cite{paper}.
In order to prove theorem \ref{theo:Gg} we must link these geometric considerations to the gradients of the
extended sections  $\tilde{f}_m$. This is the contents of the following lemma.

\begin{lm}\label{lem:comparison}
Let $x \in \Perf \setminus Z_g \cong \Vor \setminus \sigma^{-1}(Z_g)$
and let $(X_x,\Theta_x)$ be the associated polarized stable semi-abelic variety,
with an involution $j_x$ on $X_x$, and the divisor $\Theta_x$ preserved by $j_x$.
Assume the following holds: for every fixed point $P_0 \in \Theta_x$ of $j_x$,
there is a component $(X^0_x,\Theta^0_x)$ of $(X_x,\Theta_x)$ whose normalization $(\tilde X^0_x,\tilde \Theta^0_x)$ is smooth at $x$,
such that the gradient of the pullback $\tilde T^0_x$ of the equation $T_x$ of
$\Theta_x$ to $\tilde X^0_x$ does not vanish:
${\rm grad}\, \tilde T^0_x(P_0) \neq 0$. Then $P_0 \not\in \overline{G^{(g)}}$.
\end{lm}
\begin{proof}
We consider the pullback of the universal family ${\mathcal X}_g^{\op {Vor}} \to \Vor$ to $\Vor(4,8)$ and call this
${\mathcal Y}$.
Over ${\mathcal A}_g(4,8)$ the $2$-torsion points on ${\mathcal Y}$ form $2^{2g}$ disjoint
sections which extend to (no longer
disjoint) sections $Z_m$ over $\Vor(4,8)$. By assumption $P_0$ lies in (at least) one section $Z_m$ parameterizing
odd $2$-torsion over ${\mathcal A}_g(4,8)$. Our aim is to show that under
the assumptions of the lemma $\tilde f_m (P_0) \neq 0$.

Since $\tilde X^0_x$ is smooth we can choose a set of $g$ local holomorphic functions on the
component of $X^0_x$, say $u^0_1, \ldots, u^0_g$, such that their pullback to
$\tilde X^0_x$ is a regular set of local parameters.
We extend these to local holomorphic functions $u_1, \ldots, u_g$ in a neighborhood of $P_0$ on ${\mathcal Y}$.
Since $u^0_1, \ldots, u^0_g$ is a set of parameters on $\tilde X^0_x$ near $P_0$ it follows that the restriction of
$u=(u_1, \ldots, u_g)$ to smooth fibers of ${\mathcal Y}$ form a regular set of parameters at points $P \in Z_m$ in
a neighborhood of $P_0$.

We can extend the equation $T_x$ of the theta divisor on $X_x$ to an equation of the
pullback of the universal theta divisor on ${\mathcal Y}$ in nearby fibers. We denote this extension by $T$.
Then ${\rm grad}_u T(P)$ and $\tilde f_m$ are proportional after a suitable local trivialization of the
Hodge bundle $\EE$. This is true for smooth fibers of ${\mathcal Y}$, but also holds for torus rank $1$ degenerations.
This can be easily seen by computing the Fourier expansion of the sections $\tilde{f}_m$
with respect to the coordinate $t=\ee(\tau_{1,1})$, see also \cite{paper} and
\cite[Sec.~5]{grsmconjectures}.
Hence there is a nowhere vanishing function $h$ in a neighborhood of the point $P_0$ in $Z_m$
such that ${\rm grad}_u T(P)= h \cdot \tilde f_m$. Since this holds on Mumford's partial
compactification, i.e. in codimension 2 in $(\Vor)^0$, and since $\Vor(4,8)$ and hence $Z_m$ is smooth near this point we can
extend $h$ as well as
the above equality to a neighborhood of $P_0$ in $Z_m$. Since $h$ does not have zeroes in codimension $2$ it is nowhere
vanishing.
Now the claim follows from the assumption that ${\rm grad}~ \tilde T^0_x(P_0) \neq 0$.
\end{proof}

\section{Theta gradients vanishing on boundary divisors} \label{sec:thetagradients}
In \cite{grsmconjectures} the behavior of the gradients of the theta function near the boundary of the partial compactification  $\calA_g\sqcup\beta_1^0$ was studied. This was obtained by using the Fourier-Jacobi expansion of the theta function. However, this analytic method would not work for going deeper into the boundary, and in this section of the text we present the geometric approach to describing the geometry of polarized semi-abelic varieties. While the case of torus rank one goes back to Mumford \cite{mumforddimag}, and the geometry of the semi-abelic variety and its polarization is known, the method we explain in detail here will allow us to deal with much more complicated degenerations later on. Still, even for the case of torus rank one, the description of the involution and its fixed points seem to be new.

The singular varieties parameterized by
$\beta_1^0$ are $\PP^1$ bundles over a $(g-1)$-dimensional abelian variety $B$
with the $0$-section and the $\infty$-section glued with a shift. More precisely, such a
singular variety is determined by a pair $(B,b)$ where $B\in \calA_{g-1}$ and $b \in B$.
Identifying $B$ via the principal polarization with its dual variety we can interpret
$b$ as a degree $0$ line bundle $\calO_{B}(b)$ on $B$ and the $\PP^1$ bundle in
question is $\PP(\calO_{B}\oplus\calO_{B}(b))$.
The $0$ and $\infty$ sections of this bundle correspond to $\PP(0\oplus\calO_{B}(b))$
and $\PP(\calO_{B}\oplus 0)$ respectively (we use the geometric notation for the projective space).

The singular variety $X=X(B,b)$ is obtained from this $\PP^1$ bundle by identifying the $0$-section and
the $\infty$-section with a shift, to be precise $(z',0) \sim (z'-b,\infty)$. Note that $\pm b$
determine isomorphic singular varieties.
The semi-abelian variety acting on this is the
$\CC^*$ bundle which is the $\PP^1$ bundle with the sections at $0, \infty$ removed, which is nothing else
but the semi-abelian variety defined by the extension $b \in \operatorname{Ext}^1(B,\CC^*) \cong \Pic^0(B)$.

We have to understand the theta divisor on
this singular variety.
The idea, which we also want to apply to other degenerations, is to obtain an equation of
the theta divisor on the semi-abelic variety from geometric principles, identify the limits of
two-torsion points, and then
study the vanishing properties of the gradient of the theta function at these points (which is to say
the vanishing orders of $\tilde f_m$).
This simple torus rank one case will showcase the approach that we will take in dealing with semi-abelic varieties of higher torus rank.

In the case of torus rank one degenerations it was Mumford who first wrote down such a theta divisor $\Theta$
by studying
limits of the theta function.
This approach was further advanced in \cite{grle}. For the genus $2$ case also see \cite{hukawebook}.
We consider a one-parameter family of principally polarized abelian varieties degenerating to a
semi-abelic variety of torus rank one --- analytically this is the family for which one can compute the standard Fourier-Jacobi expansion, as in \cite{paper}.
An analytic expression for the theta function on the semi-abelic variety in this family is known.
Following Mumford, one considers a translate of the Riemann theta function by shifting the argument $z$ by
$(-w/2,0, \ldots, 0)$, i.e.~by replacing $z=(z_1,z')$ by $(z_1 - w/2,z')$ (where
$w=\tau_{1,1}$). Expanding this function in Fourier series we obtain
$$
\theta\left(\begin{pmatrix} \omega& ^{t}b \\ b &\tau' \end{pmatrix},\begin{pmatrix} z_1 - \omega/2\\ z'
\end{pmatrix}\right)
$$
$$
\begin{aligned}
&=\sum\limits_{N\in\ZZ,n'\in\ZZ^{g-1}}\ee\left(\frac12
\left(\begin{matrix}N\\ n'\end{matrix}\right)^t\left(\begin{matrix} \omega& ^{t}b \\ b &\tau' \end{matrix}
\right)\left(\begin{matrix} N\\ n'\end{matrix}\right)+\left(\begin{matrix}N\\ n'\end{matrix}\right)^t\left(\begin{matrix} z_1 - \omega/2\\ z'\end{matrix}\right)\right)
\\
&=\sum\limits_{N\in\ZZ} q^{ \frac{1}{2}N(N-1)}\ee(Nz_1) \left( \sum\limits_{n'\in\ZZ^{g-1}} \ee(\frac12 n'^t\tau'n' + n'^{t}(z'+Nb)) \right),
\end{aligned}
$$
where as before $q=\ee(w)$.
As $\omega \to i\infty$ we have $q \to 0$, so that only the terms with $N=0$ and $N=1$ survive, and thus as a limit we obtain the semi-abelic theta function
\begin{equation}\label{thetamumford}
  T(z',x):=\theta(\tau',z') + x \theta(\tau',z'+b),
\end{equation}
where $x:=\ee(z_1)$ can be viewed as the fiber coordinate of the $\CC^*$ or $\PP^1$ bundle,
and the gluing is $(z',0)\sim(z'-b,\infty)$. This gives an analytic equation for $\Theta$.

To simplify notation, from now on we will drop the $'$ and denote $z'$ and $\tau'$ simply by
$z$ and $\tau$ when no  confusion is possible; this convention will also be used in the following sections.

We now come to the delicate point of determining the involution on the
semi-abelic pair $(X,\Theta)$, and determining the limits of two-torsion points on it as the
fixed points of this involution --- to the best of our knowledge these questions have not been previously treated in the literature. It follows from general theory that we
have an involution $j$ on the polarized semi-abelic variety which is
compatible with the inverse map on the semi-abelian (open) variety, i.e. the action of the
semi-abelian variety on $(X,\Theta)$ is equivariant with respect to the involution $j$ on $X$ and the inverse
map of the semi-abelian variety viewed as a group scheme. We want to write this down explicitly.
However, one needs to be careful in writing down the involution in $z,x$
coordinates. The issue here is that while we are interested in the
involution on a semi-abelic variety, and its fixed points there,
analytically the computation happens on the universal cover.

Indeed, we are interested in the involution on the $\PP^1$ bundle over the
abelian variety $B$, while formula (\ref{thetamumford}) is in fact an analytic
expression for the lifting of the theta function from the semi-abelian variety
to the universal cover $\CC^*\times \CC^{g-1}$. To describe the involution explicitly we first
describe a uniformization of the smooth part of $X$ (which can be identified with the semi-abelian variety
corresponding to the point $b$). Let
$B=\CC^{g-1}/(\ZZ^{g-1}+\ZZ^{g-1}\tau)$. The semi-abelian variety corresponding
to the line bundle $b\in B^\vee$ is the quotient of the  trivial $\CC^*$ bundle
$\CC^{g-1}\times \CC^*$ over $\CC^{g-1}$ by the group
$\ZZ^{2(g-1)} \cong\ZZ^{g-1}+\ZZ^{g-1}\tau$ acting via
\begin{equation}\label{transform}
 (n,m):(z,x)\mapsto (z+n\tau+m,\ee(-n^tb)x),
\end{equation}
where $(n,m) \in \ZZ^{2(g-1)}$ (see for example \cite[Sec.~II.3]{igusabook} or \cite[Sec.~I.2]{bila})

To determine the action of the involution $j$ on the semi-abelic variety, we note that
the involution must lift to an involution on $\CC^{g-1}\times\CC^*$, inducing an
involution on the base, $z\mapsto a-z$ for some $a\in B$, and mapping fibers to
fibers. If we complete the $\CC^*$
bundle to a $\PP^1$ bundle, the involution must interchange the $0$ and $\infty$
sections, and thus must be given by $(z,x)\mapsto (a-z,c(z)x^{-1})$ on the fiber over $z\in\CC^{g-1}$
for some $c:\CC^{g-1}\to\CC^*$. For this involution to be preserved by the action (\ref{transform})
we must then have $c(z+n\tau+m)=c(z)$ for any $n,m\in\ZZ^g$, and thus
$c$ must be  a holomorphic function $c:B\to\CC^*$, and thus constant.

We now determine the constants $a\in B$ and $c\in \CC^*$ by requiring the involution
to preserve the semi-abelic theta divisor. Indeed, we have
$$
  T(j(z,x)) = T(a-z,cx^{-1})=\theta(\tau,a-z)+cx^{-1}\theta(\tau,a-z+b)=
$$
$$
  \theta(\tau,z-a)+cx^{-1}\theta(\tau,z-a-b)=cx^{-1}\left(\theta(\tau,z-a-b))
  +c^{-1}x\theta(\tau,z-a)\right)
$$
and for this to match $T(z,x)$ we must then have $a=-b$ and $c=1$, as expected (recall
that in Mumford's construction the origin was shifted by $b/2$).

We now determine the limit of two-torsion points on the degenerating family of
abelian varieties --- these are the fixed points of the involution $j$ on the
semi-abelic variety. These points can either lie in the smooth part of $X$ or in its singular
locus (which is isomorphic to the abelian variety $B$). We start with the smooth locus.
Note that a fixed point of the
involution $j$ on the $\CC^*$ bundle over $B$ does not necessarily lift to a
fixed point of the involution on $\CC^{g-1}\times \CC^*$; rather it lifts to a
point such that its image under the involution is identified with the point itself
under the action (\ref{transform}).
Indeed, suppose some $(z,x)\in\CC^{g-1}\times\CC^*$ descends to a fixed point
of $j$ on the $\CC^*$ bundle over $B$. This means that there exist
$\ep,\de\in\ZZ^{g-1}$ such that
$$
 j(z,x)=(-b-z,x^{-1})=(z+\tau\ep+\de,\ee(-\ep^tb)x).
$$
This implies that $z=-(\tau\ep+\de+b)/2$ (as a point in $\CC^{g-1}$, not on $B$!)
and that $x=\pm\ee(-\ep^tb/2)$. Note that since $\ep$ and $\de$ differing by even integer vectors give the
same point in $B$, it is enough to consider $\ep,
\de\in (\ZZ/2\ZZ)^{g-1}$, and thus the total number of
fixed points of $j$ we obtain in this way is $2\cdot 2^{2g-2}$. Each of them arises as
a limit of a single family of two-torsion points on smooth abelian varieties (geometrically
this means that the multisection defined by these $2$-torsion points on the universal family extends locally to sections
over the boundary and is unramified there).

The fixed points of the involution on the glued $0$ and $\infty$ sections of the
semi-abelic variety can be obtained from the equation
$$
 j(z,0)=(-b-z,\infty)\sim (-z,0)=(z+\tau\ep+\de,0),
$$
where $\sim$ denotes the glued points on the $0$ and $\infty$ sections. This implies
that $z=-(\tau\ep+\de)/2$.
In this way we have found $2^{2g-2}$ fixed points of the involution on the singular locus of the semi-abelic variety.
Each of these
points is the limit of two $2$-torsion points on a smooth abelian variety
(the extended multisection of $2$-torsion points is ramified here) and hence has to be counted with multiplicity $2$.
Together with the $2^{2g-1}$ fixed points on the smooth part we indeed obtain $2^{2g}$ points.

\smallskip
Note that for the simple example of an elliptic curve degenerating as $\tau\to i\infty$,
in the limit the one fixed point of the involution on the glued locus is a limit of both two-torsion points
$\tau/2$ and $(1+\tau)/2$. Thus in a family of
abelian varieties degenerating to a semi-abelic variety it does not make sense to say that some
fixed points of the involution on the semi-abelic variety are even and some are odd. However, it
makes sense to ask which fixed points are limits of odd points (and possibly also of some even ones) --- these are the fixed points of $j$ that generically lie on the semi-abelic theta
divisor. We will now determine these points and then compute
the gradients (with respect to $z$ and $x$) of the semi-abelic theta function $T(z,x)$
at these fixed points lying on the semi-abelic theta divisors, and verify that these
gradients generically do not vanish. As we have explained, this
gives us information on the vanishing locus of the section $\tilde f_m$.
The gradient of the torus rank 1 semi-abelic theta function at a
generic point was computed in \cite{civdg2}, where the singularities of the theta divisor
were studied --- while we are only interested in the singularities at the fixed points of the involution $j$.

For points in the smooth part of the semi-abelic variety we denote $m:=(\tau\ep+\de)/2$
the two-torsion point on $B$. Using the fact that the theta function is even we compute
$$\begin{aligned}
 T(-(\tau\ep+\de+b)/2,&\pm\ee(-\ep^tb/2))\\
 &=\theta(\tau,-m-b/2)\pm\ee(-\ep^tb/2)\theta(\tau,-m+b/2)\\
 &=\theta(\tau,m+b/2)\pm\ee(-\ep^tb/2)\theta(\tau,-m+b/2).
\end{aligned}
$$
Now we note that $2m=\tau\ep+\de$ and so by quasi-periodicity of the theta function
we get
$$
\begin{aligned}
 \theta(\tau,m+b/2)&=\theta(\tau,-m+b/2+\tau\ep+\de)\\
&=\ee(-\ep^t\tau \ep/2-\ep^t(-\tau \ep+\de+b)/2 )\theta(\tau,-m+b/2)\\
 &=(-1)^{\ep^t\de}\ee(-\ep^tb/2)\theta(\tau,-m+b/2).
\end{aligned}
$$
Substituting this into the above expression finally yields
$$
 T(-m-b/2,\pm\ee(-\ep^tb/2))=\ee(-\ep^tb/2)\theta(\tau,-m+b/2)((-1)^{\ep^t\de}\pm 1),
$$
and thus the semi-abelic theta function vanishes at this fixed point of $j$ generically
if $m\in B[2]_{\rm odd}$, and the plus sign is chosen, or if $m\in B[2]_{\rm even}$,
and the minus sign is chosen.

We now compute the gradient of the semi-abelic theta function at these points.
Computing these gradients is the same as evaluating the sections $\tilde{f}_m$,
as we saw in lemma \ref{lem:comparison}.

Indeed, we have
$$
 \frac{\partial}{\partial x}T(z,x)=\theta(\tau,z+b);\qquad
 \frac{\partial}{\partial z}T(z,x)=\frac{\partial\theta(\tau,z)}{\partial z}
 +x\frac{\partial\theta(\tau,z+b)}{\partial z},
$$
and plugging in an odd two-torsion point gives these derivatives as
\begin{equation}\label{equ:partial1}
 \theta(\tau,-m+b/2);\qquad 2\ee(-\ep^tb/2)\frac{\partial\theta}{\partial z}(\tau,-m+b/2)
\end{equation}
using quasi-periodicity as above. For this gradient to vanish is then equivalent
to the point $-m+b/2$ being a singular point on the theta divisor of $B$, and this
condition is known to define a codimension $g$ locus within the universal family $\calX_{g-1}$ of $(g-1)$-dimensional abelian varieties.

To finish the discussion of the torus rank one case, we need to determine which fixed
points of the involution on the singular part of the semi-abelic variety lie on
the theta divisors and compute the gradient of theta there. Indeed, we have
$$
 T(-m,0)=\theta(\tau,-m),
$$
which vanishes identically if and only if $m$ is an odd two-torsion point. The
gradient of the semi-abelic theta function at such a point is
\begin{equation}\label{equ:partial2}
 \theta(\tau,-m+b); \qquad \frac{\partial\theta}{\partial z}(\tau,-m).
\end{equation}
The condition for the derivatives to vanish is exactly to say that $\tau\in I^{(g-1)}$.
Then the condition for $\theta(\tau,-m+b)$ to vanish is equivalent to the corresponding
theta function with characteristic vanishing at $b$, which is clearly an independent
vanishing condition, so that in this case we arrive at the following

\begin{prop}\label{codim1Ok}
If theorem \ref{theo:Ig} holds in genus $g-1$, then in genus $g$ the zero locus of each $\tilde f_m$ does not have any irreducible components contained in $\beta_1^0$. It thus follows that we have
$$
 \overline{I^{(g)}}\cap(\calA_g\sqcup\beta_1^0)= \overline{G^{(g)}}\cap(\calA_g\sqcup\beta_1^0),
$$
i.e.~that the locus $\overline{G^{(g)}}$ does not pick up any ``extra components''
for semi-abelic varieties of torus rank 1.
\end{prop}
\begin{rem}
In the following sections of the text we will study the gradients of the semi-abelic theta functions for various types of principally polarized semi-abelic varieties. The conclusion in each case will be the same: that if theorems \ref{theo:Ig} and \ref{theo:Gg} hold in lower genus, then in genus $g$ the loci $\overline{G^{(g)}}$ and $\overline{I^{(g)}}$, when intersected with the stratum of semi-abelic varieties of the type considered, coincide and are both of codimension $g$. Altogether, this will eventually yield the proofs of theorems \ref{theo:Ig} and \ref{theo:Gg}.
\end{rem}

\section{Theta gradients on the standard degenerations} \label{sec:standardcomp}
The easiest example of a semi-abelic variety is when $\RR^n$ is divided into standard cubes with
$\ZZ^n$ acting by translation along the coordinate axes.
The following figure (see Figure \ref{Figure1}) depicts this for $n=2$.

\begin{figure}[ht!]
\centering
\includegraphics[scale=0.8]{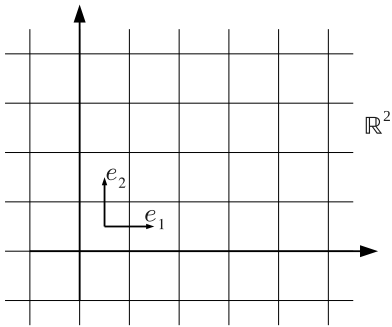}
\caption{The standard dicing of the plane with translations}\label{Figure1}
\end{figure}

The normalization of the semi-abelic variety is then a $(\PP^1)^n$ bundle over an abelian variety $B$ of dimension $g-n$.
Following \cite[Sec.~16ff]{civdg2} we call such a semi-abelic variety {\it standard}.
For $\Perf$ (but not for $\Vor$!) the locus
parameterizing standard semi-abelic varieties is open and dense in $\beta_n$, see \cite{shepherdbarron}.
These degenerations of abelian varieties have already been studied by different authors, see
\cite[Sec.~4ff.]{dhs}, \cite[Sec.~16ff.]{civdg2} (the conventions in these texts are slightly different, resulting in
different signs; we make a choice compatible with Mumford's description in the torus rank one case).
We shall recall the
description of these semi-abelic varieties and
of their generalized theta divisors, and compute the gradient of the theta divisor, which
allows us to compute the extension of $\tilde f_m$ to this
locus.

In terms of period matrices
these degenerations correspond to the degenerations of the period matrix where the first $n$
diagonal entries tend to $i\infty$, and the others stay bounded:
$$
 \tau=\begin{pmatrix}\tau_{11}&\ldots&\tau_{1k}&\ldots&\tau_{1n}&b_1^t\\
 \vdots&\ddots&\ldots&\ldots&\ldots&\ldots\\
 \tau_{j1}&\vdots&\tau_{jj}&\ldots&\tau_{jn}&b_j^t\\
 \vdots&\vdots&\vdots&\ddots&\ldots&\ldots\\
 \tau_{n1}&\vdots&\tau_{nk}&\vdots&\tau_{nn}&b_n^t\\
 b_1&\vdots&b_k&\vdots&b_n&\tau'\end{pmatrix}
$$
with $\tau_{jj}\to i\infty$ for $1\le j\le n$, $\tau\in\calH_{g-n}$ being the period
matrix of $B$, and $b_1,\ldots,b_n\in B$.
The semi-abelic variety is then obtained from the $(\PP^1)^n$ bundle
over $B\in\calA_{g-n}$ (the non-homogenous fiber coordinates on which we denote $x_1,\ldots,x_n\in\PP^1$) by gluing
$$
 (z,x_1,\ldots,x_{j-1},0,x_{j+1},\ldots,x_n)\sim
$$
$$
\sim (z-b_j,t_{j,1}^{-1}x_1,\ldots,t_{j,j-1}^{-1}x_{j-1},\infty,
 t_{j,j+1}^{-1}x_{j+1},\ldots,t_{j,n}x_n)
$$
for all $1\le j\le n$, where the parameters $t$ are given by $t_{j,k}=\ee(\tau_{jk})$.
Thinking of the toric picture of $(\PP^1)^n$ as a
hypercube, this means that each of the $n$ pairs of parallel $(n-1)$-dimensional faces is identified with a shift.

The semi-abelic theta divisor was then computed in \cite{dhs,civdg2} to be given by
$$
 T(z,x_1,\ldots,x_n)=\sum\limits_{\mu_1,\ldots,\mu_n\in\lbrace 0,1\rbrace}\prod\limits_{j=1}^n x_j^{\mu_j}
 \prod\limits_{1\le j<k\le n}t_{j,k}^{\mu_j\mu_k} \theta\left(\tau,z+\sum\limits_{j=1}^n \mu_jb_j\right).
$$
By tedious computations generalizing the case of $n=1$, the involution on this semi-abelic variety
can be computed to be given by
$$
 j(z,\ldots,x_j,\ldots)=(-\sum b_j-z,\ldots,x_j^{-1}\prod\limits_{k\ne j} t_{j,k}^{-1},\ldots).
$$
Instead of going through the derivation, let us check that this works, i.e.~that this involution
maps glued points to glued points, and preserves the semi-abelic theta divisor. Indeed, the point
$$
 j(z-b_j,t_{j,1}^{-1}x_1,\ldots,\infty,\ldots,)=
$$
$$
=(-b_1-\ldots-b_n+b_j-z,(x_1t_{j,1}^{-1})^{-1}\prod\limits_{k\ne 1} t_{1,k}^{-1},
 \ldots,0,\ldots)
$$
is glued to
$$
 (-b_1-\ldots-b_n+b_j-z-b_j,t_{j,1}^{-1}(x_1t_{j,1}^{-1})^{-1}\prod\limits_{k\ne 1} t_{1,k}^{-1},\ldots,\infty,\ldots),
$$
which is equal to $j(z,x_1,\ldots,0,\ldots)$, and thus the involution preserves the gluing. We now
check that the semi-abelic theta divisor is invariant under $j$. Indeed, we have
$$
 T(j(z,x_1,\ldots))=
$$
$$
\sum\limits_{\mu_1,\ldots,\mu_n\in\lbrace 0,1\rbrace} \prod\limits_{j=1}^n
 \left(x_j^{-1}\prod\limits_{k\ne j}t_{j,k}^{-1}\right)^{\mu_j}\prod\limits_{1\le j<k\le n}
 t_{j,k}^{\mu_j\mu_k}\theta(-z+\sum (\mu_j-1)b_j)
$$
$$
 =\prod\limits_{j=1}^n x_j^{-1}\prod\limits_{1\le j<k\le n}t_{j,k}^{-1}
 \sum\limits_{\mu_1,\ldots,\mu_n}\prod\limits_{j=1}^n x_j^{1-\mu_j} \times
$$
$$
\ \  \ \ \ \prod\limits_{1\le j<k\le n}
 t_{j,k}^{\mu_j\mu_k-\mu_j-\mu_k+1}\theta(z+\sum(1-\mu_j)b_j)
$$
which upon changing the variables of summation to $1-\mu_j$ becomes simply the
expression for
$$
 \prod\limits_{j=1}^n x_j^{-1}\prod\limits_{1\le j<k\le n}t_{j,k}^{-1}\ T(z,x_1,\ldots,x_n).
$$
To determine the fixed points of $j$, we again need to be careful to distinguish working on
$(\CC^*)^n\times\CC^{g-n}$ and on the bundle over $B$. The same arguments as for the $n=1$ case
show that the fixed points of the involution $j$ are of the form
$$
 \left(-\frac{\tau\ep+\de+\sum \mu_i b_i}{2},\ldots,\pm\mu_j\ee(-\ep^tb_i/2)\prod\limits_{k\ne j}t_{j,k}^{-1/2},\ldots\right),
$$
where $\mu_1,\ldots,\mu_n\in\lbrace 0,1\rbrace$.

Similarly to the case of torus rank one, one can then determine which of the fixed
points of $j$ are ``odd'', i.e.~which of them generically lie on the semi-abelic theta divisor,
and then one verifies that the gradient of $T$ (i.e.~the set of derivatives with respect to $z,x_i$)
impose $g$ independent conditions (assuming that theorems \ref{theo:Ig} and \ref{theo:Gg} hold for $g-n$).
The computation is straightforward, but tedious, and we skip it; the $n=2$ case easily follows
from the computations of gradients done in \cite{civdg2}.

\section{The non-standard semi-abelic variety of torus rank $2$}

\begin{figure}[ht!]
\includegraphics[scale=0.8]{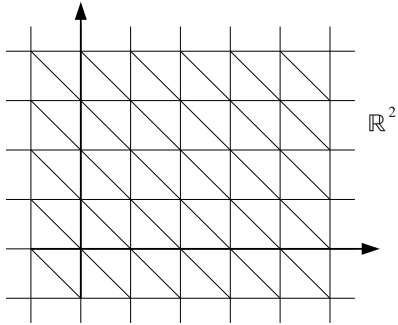}
\caption{The non-standard dicing of the plane} \label{Figure1a}
\end{figure}

For torus rank $2$ there is only one other type of degeneration apart from the standard degeneration discussed in the
previous section.
In this case the
decomposition of $\RR^2$ is obtained by subdividing each unit square into two triangles (see Figure \ref{Figure1a}).
Hence the normalization of the semi-abelic variety consists of two $\PP^2$ fibrations
over an abelian variety $B\in\calA_{g-2}$.
The locus of the non-standard rank $2$ degenerations in $\Perf$
is the so-called singular locus $\Delta$ in $\beta_2^0$, in notation of \cite{vdgeercycles} and \cite{ergrhu2}.
It has codimension one in $\beta_2^0$, and in terms of the above discussion corresponds to the
case of $t_{1,2}=\infty$, so that the previous computations break down.
The Fourier-Jacobi expansion in this case is very complicated to deal with.
We shall, therefore, apply our general geometric approach described in section \ref{sec:outline}.

The description of the geometry in this case can also be obtained from the general construction in \cite{alna,alexeev}, and the
entire picture
of the semi-abelic variety and the semi-abelic theta function were computed and written out
in detail in \cite[Sec.~17]{civdg2}, where it was also considered as a $t_{1,2}\to \infty$ degeneration of the
standard compactification. Our reason for giving now a completely elementary derivation for this semi-abelic theta divisor is
to showcase all the techniques that we will then use for the much harder and not previously dealt with, non-standard torus rank $3$
degenerations. Our sign conventions are slightly different
from those of \cite{civdg2} --- we chose the conventions to be consistent throughout this paper.

\begin{figure}[ht!]
\centering
\includegraphics[scale=0.8]{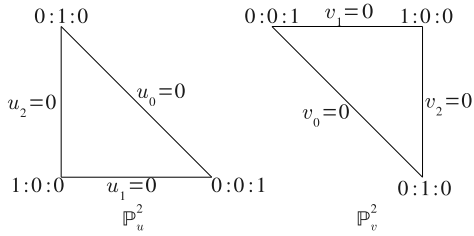}
\caption{The gluings for the non-standard torus rank $2$ case}\label{Figure2}
\end{figure}

The normalization of the semi-abelic variety in this case consists of two copies of a $\PP^2$ bundle over
some $B\in\calA_{g-2}$. We will denote by $z\in B$ the coordinates on the bases of these bundles,
and by $(u_0:u_1:u_2)$ and $(v_0:v_1:v_2)$ the homogeneous coordinates on the fibers $\PP^2_u$ and $\PP^2_v$ of these bundles.
In the purely toric case (i.e.~genus $2$)
the two copies of the $\PP^2$ bundle are glued along the $3$ coordinate $\PP^1$'s.
We look at the corresponding toric picture of polytopes: two triangles glued into a square, with the opposite edges further identified, and
with the involution
being the symmetry with respect to the center of the square --- see Figure \ref{Figure2}.
Then the edges of the $\PP^2$ correspond to
the coordinate lines $\PP^1$ obtained by setting one homogeneous coordinate to zero, and the vertices to the axis points
on the coordinate lines.

\begin{figure}[ht!]
\centering
\includegraphics[scale=0.8]{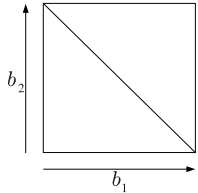}
\caption{Shifts in the non-standard torus rank $2$ case} \label{Figure3}
\end{figure}

In general (i.e. in the presence of a non-trivial abelian base),
fibers over the points in the base differing by shifts (i.e.~by adding some point on $B$) are glued as indicated in Figure \ref{Figure3}.
When identifying the $u_0=0$
and $v_0=0$ coordinate $\PP^1$'s, there is no shift on $B$ (this being an interior face of the
decomposition of the cube), while the shifts for the other two gluings are $b_1$ and $b_2$, the same as for the case of the standard compactification $(\PP^1)^2$.

Now to understand the gluing parameters for the fiber coordinates, in $\PP^2_u$ and $\PP^2_v$'s, we have a choice. The simplest choice here is to require the involution to map the coordinates to the
coordinates (which corresponds to the labeling of the coordinate lines in the two $\PP^2$'s as in the figure). The involution is then given in coordinates as
$$
 j:(z,(p:q:r)_u)\mapsto (-b_1-b_2-z,(p:q:r)_v)
$$
where the subscript indicates on which $\PP^2$ the point is taken.

To see what the gluings of the three coordinate $\PP^1$'s on $\PP^2_u$ to those on $\PP^2_v$ are in these
coordinates, we look at the vertices.
For example, we see from figure 3 that the gluing identifying the diagonal sides of the triangles glues
the vertex corresponding to $(0:0:1)_u$ to $(0:1:0)_v$, and glues $(0:1:0)_u$ to $(0:0:1)_v$.
It thus follows that the $\PP^1$ given by $u_0=0$ on $\PP^2_u$ is glued  to the $\PP^1$
given by $v_0=0$ on $\PP^2_v$ with an involution, and thus the gluing is
given by $(z,(0:p:q))_u\to(z,(0:\lambda q:p))_v$, for some parameter $\lambda \in\CC^*$. Since we still have
the freedom of choosing the coordinates on $\PP^2_u$, we can rescale the coordinate $u_1$ so that this
parameter is made equal to $1$
(the situation will be much more complicated for the case of torus rank $3$ that we will need to consider
in what follows).
Similarly by rescaling the coordinate $u_2$ (and simultaneously also $v_2$ which is identified with $u_2$
by the involution $j$) we can make the second gluing to be $(z,(p:0:q)_u)\to (z-b_1,(q:0:p)_v)$,
whereas for the third gluing we can, a priori, only say that it is of the form
$(z,(p:q:0))_u)\to (z-b_2,(c q:p:0)_v)$ for some gluing parameter $c \in \CC^*$.
We shall see below that the existence of a symmetric principal polarization implies $c = 1$.
This description was obtained in \cite[Sec.~17]{civdg2} (with which we have a different sign convention)
by studying the appropriate identifications of the summands of rank three vector bundles over $B$ projectivized
to the $\PP^2$ bundles.

The theta divisor on the normalization of the semi-abelic variety has the following form. On the section given by
the coordinate points of the $\PP^2$ bundles it is isomorphic to $\Theta_B$, whereas its restriction to each
$\PP^2$ in a fiber is a section of $\calO(1)$, as follows from the general theory. It can also be seen in a more elementary way
by looking at the split situation ($b_1=b_2=0$). Then this is the only choice of a class of an ample line bundle
with correct top self-intersection number.

Thus if we denote by $T_u$ and $T_v$ the classes of the theta divisors on the two copies of the $\PP^2$ bundle (with coordinates on the fibers being $u$ and $v$, respectively),
we must have (where we dropped the argument $\tau_B$ of the theta function everywhere)
$$
 T_u(z_1,(u_0:u_1:u_2))=\lambda_0u_0\theta(z_1+\alpha_0)+\lambda_1u_1\theta(z_1+\alpha_1)+\lambda_2u_2\theta(z_1+\alpha_2);
$$
$$
 T_v(z_2,(v_0:v_1:v_2))=\mu_0v_0\theta(z_2+\beta_0)+\mu_1v_1\theta(z_2+\beta_1)+\mu_2v_2\theta(z_2+\beta_2)
$$
for some points $\alpha_i,\beta_i\in B$ (in \cite{civdg2} it is explained why these shifts are in
fact known a priori, but here we give a direct way of computing them from the gluings), and
for some coefficients $\lambda_i,\mu_i\in\CC$.
The involution $j$ must interchange $T_u$ and $T_v$, which is equivalent to requiring
$$
 \theta(z+\alpha_i)\sim\theta(-b_1-b_2-z+\beta_i)=\theta(z+b_1+b_2-\beta_i),
$$
where by abuse of notation $\sim$ denotes equality up to a non-zero constant (we will usually
get such identities from the gluings, also denoted $\sim$, so the notation makes sense).

Thus for the shifts we have $\alpha_i=b_1+b_2-\beta_i$ and for the coefficients we have $\lambda_i=\mu_i$.
We now write down the conditions for $T_u$ and $T_v $ to glue, i.e.~for the
divisors to coincide on the glued coordinate $\PP^1$'s. Indeed, we must have
$$
 T_u(z,(0:p:q)) \sim T_v(z,(0:q:p)),
$$
which means that we must have
$$
 \lambda_1 p\theta(z+\alpha_1)+\lambda_2 q\theta(z+\alpha_2) \sim \lambda_1q\theta(z+b_1+b_2-\alpha_1)+\lambda_2p\theta(z+b_1+b_2-\alpha_2)
$$
for all $z\in B$, and all $u,v\in\CC$. This is equivalent simply to
$$
 \alpha_1=b_1+b_2-\alpha_2;\qquad (\lambda_1:\lambda_2)=(\lambda_2:\lambda_1).
$$
Similarly from
$$\begin{aligned}
 T_u(z,(p:0:q))&=\lambda_0 p\theta(z+\alpha_0)+\lambda_2q\theta(z+\alpha_2)\\
 \sim T_v(z-b_1,(q:0:p))&=\lambda_0q\theta(z+b_2-\alpha_0)+\lambda_2p\theta(z+b_2-\alpha_2)
 \end{aligned}
$$
we get
$$
 \alpha_0=b_2-\alpha_2;\qquad (\lambda_0:\lambda_2)=(\lambda_2:\lambda_0),
$$
and finally from
$$\begin{aligned}
 T_u(z,(p:q:0))&=\lambda_0p\theta(z+\alpha_0)+\lambda_1q\theta(z+\alpha_1)\\
 \sim T_v(z-b_2,(c q:p:0))&=c \lambda_0q\theta(z+b_1-\alpha_0)+\lambda_1p\theta(z+b_1-\alpha_1)
 \end{aligned}
$$
we get
$$
 \alpha_0=b_1-\alpha_1; \qquad (\lambda_0:\lambda_1)=(c\lambda_0:\lambda_1).
$$
Notice that combining this proportionality with the previous ones gives $\lambda_0=\lambda_1=\lambda_2$
and $c=1$.
Solving the three equations for the shifts we finally obtain
$$\begin{aligned}
 T_u(z,(u_0:u_1:u_2))&=u_0\theta(z)+u_1\theta(z+b_1)+u_2\theta(z+b_2);\\
 T_v(z,(v_0:v_1:v_2))&=v_0\theta(z+b_1+b_2)+v_1\theta(z+b_2)+v_2\theta(z+b_1).
 \end{aligned}
$$
We now compute the fixed points of the involution $j$; since it interchanges the two components,
the fixed points must lie on the gluing locus, i.e.~on the union of coordinate $\PP^1$'s.
We first compute the fixed points on the three gluing affine lines, and then deal with the vertices of the triangle.
As in the previous computations, one delicate point is that we are working on a non-trivial  $(\CC^*)^2$ bundle over $B$.
Indeed, we shall make repeated use of the $2$-dimensional analogue of formula (\ref{transform}).

On the line $u_0=0\ne u_1u_2$ the fixed points of the involution $j$ are  given by the condition that
$$
 j((z,(0:p:q))_u)=(-b_1-b_2-z,(0:p:q))_v
$$
is identified with itself under the gluing $(z,(0:p:q))_u\sim(z,(0:q:p))_u$ (after a possible
shift by some
element of $\ZZ^{g-2}+\tau\ZZ^{g-2}$, and the corresponding factors appearing in $v$ and $u$).
This implies that these points are of the form $(m-(b_1+b_2)/2,(0:\pm \ee(\ep^t(b_1-b_2)/2):1)$
where $m=(\tau\ep+\de)/2\in B[2]$.
Similarly the fixed points of $j$ on the other two glued $\PP^1$'s (away from the vertices of the triangle,
i.e.~away from the three coordinate points in $\PP^2_s$) are of the form $(m-b_2/2,(\pm \ee(\ep^tb_2/2):0:1)_u)$
and $(m-b_1/2,(\pm \ee(\ep^tb_1/2):1:0)_u)$.

To determine which of these two torsion points generically lie on the theta divisor, we compute,
denoting $p:=(m-(b_1+b_2)/2,(0:\pm \ee(\ep^t(b_1-b_2)/2):1))_u$
$$
 T_u(p)=\pm \ee(\ep^t(b_1-b_2)/2)\theta(m+(b_1-b_2)/2)+\theta(m+(b_2-b_1)/2),
$$
which by using again the quasi-periodicity (\ref{qperiodic}) is proportional to
$$
 \pm\theta_m((b_1-b_2)/2)+\theta_m((b_2-b_1)/2),
$$
and vanishes if the signs are chosen appropriately.
Thus onthis $\CC^*$ the fixed points of the involution generically lying on the semi-abelic theta divisor are
$$
 \left\lbrace(m-(b_1+b_2)/2,(0:\ee(\ep^t(b_1-b_2)/2):1))\mid m\in B[2]_{\rm odd}\right\rbrace
$$
and
$$
\left\lbrace(m-(b_1+b_2)/2,(0:-\ee(\ep^t(b_1-b_2)/2):1))\mid m\in B[2]_{\rm even}\right\rbrace.
$$
We now compute the gradient of $T_u$ at any such point, i.e.~we need to compute the derivatives with respect to local coordinates $z\in B$, and $(u_1:u_2)\in\PP^1$. Again using (\ref{qperiodic}), up to an exponential factor we get
$$
 \partial_zT_u|_p=\pm\partial_z\theta_m((b_1-b_2)/2)+\partial_z\theta_m((b_2-b_1)/2),
$$
which for each of the two cases  gives simply $2\partial_z\theta_m((b_2-b_1)/2)$. For the other two partial derivatives we compute (again up to exponential factors)
$$
 \partial_{u_1}T_u|_p=\theta_m((b_1+b_2)/2);\qquad
 \partial_{u_2}T_u|_p=\theta_m((b_1-b_2)/2)
$$
(compare to the computation in \cite{civdg2}). For all of these partial derivatives to vanish
simultaneously, the point $(b_2-b_1)/2$ must be a singular point of the theta divisor $\Theta_B$
(which imposes $g-1$ conditions in $\calX_{g-2}$), and in addition we must
have $(b_1+b_2)/2\in\Theta_B$, which is another independent condition --- thus in this case we have a
codimension $g$ locus within $\Delta\subset\beta_2^0$ where the gradient of the equation defining the
theta divisor vanishes, as expected.
The computation for the fixed points of $j$ lying on the two other glued coordinate $\PP^1$'s is completely
analogous, with points $b_1/2$ and $b_2/2$ appearing instead of $(b_1+b_2)/2$, and we do not give it here.

The remaining case to handle is of the fixed points over the vertices of the triangles in the toric picture.
In this case there are multiple gluings to take care of, so we determine the fixed points of the involution as follows:
$$\begin{aligned}
 j((z,(p:0:0))_u)&=(-b_1-b_2-z,(p:0:0))_v\sim (-z-b_1,(0:p:0))_u\\
 &\sim (-z+b_2,(0:0:p))_u \sim (-z,(p:0:0))_u
 \end{aligned}
$$
where we used the gluings along all the coordinate $\PP^1$'s. Thus the fixed points in this case are
given by $(m,(1:0:0))$ for $m\in B[2]$. Note that each such fixed point is a limit of $4$ two-torsion points
on a smooth ppav. Since
$$
 T_u(m,(1:0:0))=\theta(m),
$$
such a fixed point of involution $j$ is odd (generically lies on $T_u$) if and only if $m$ is odd.
We then compute the gradient: indeed
$$
 \partial_zT_u|_{(m,(1:0:0))}=\partial_z\theta(m)
$$
vanishes if and only if $\tau\in I^{(g-2)}$, and furthermore the vanishing of
$$
 \partial_{u_1}T_u|_{(m,(1:0:0))}=-\theta(m+b_1)\qquad{\rm and}\qquad  \partial_{u_2}T_u|_{(m,(1:0:0))}=-\theta(m+b_2),
$$
are the conditions for $b_1$ and $b_2$ to lie on the translated theta divisor --- and thus give two more
independent conditions, as expected.

\section{The codimension 4 strata of non-standard semi-abelic varieties in $\Perf$}
In this section we deal with the cones that correspond to strata of non-standard semi-abelic varieties of codimension 4 in $\Perf$. Recall from the table that there are two such strata, both of them in $\beta_3^0$. The easier case, for which we can easily deduce the appropriate formulae, corresponds to
dicing the cube into two triangular prisms (so the corresponding toric variety is the union of two copies of $\PP^1\times\PP^2$). The other case - the most complicated that we need to handle in this paper --- is that of the dicing
into two tetrahedra and one octahedron (so the corresponding toric variety is the union of two copies
of $\PP^3$ and a toric threefold of type $F(2,2)$ --- the intersection
of two special quadrics in $\PP^5$).

\subsection{Two copies of a $\PP^1\times\PP^2$ bundle}
The simpler case is that of two triangular prisms, i.e.~that of a semi-abelic variety having
two irreducible components, the normalization of each being a fibration over some $B\in\calA_{g-3}$ with
fibers $\PP^1\times\PP^2$. As before for the case of two copies of $\PP^2$, we denote these two copies by
subscripts $u$ and $v$ respectively, choose coordinates on the second copy to be the image under the involution $j$ of the coordinates
on the first copy, and finally denote by $x$ the (non-homogenous) coordinate on $\PP^1_u$. As for the case
of two copies of $\PP^2$, which, as explained in \cite{civdg2}, can be obtained by degenerating the
standard $\PP^1\times\PP^1$ compactification when
letting the gluing parameter $t_{1,2}$ approach zero, this compactification can be obtained from
the standard $(\PP^1)^3$ compactification by letting $t_{1,2}$, now one of the three gluing parameters,
approach zero. We then get
$$
 T_u(x,(u_0:u_1:u_2))=u_0\theta(z)+u_1\theta(z+b_1)+u_2\theta(z+b_2)\hskip25mm
$$
\begin{equation}\label{p1p2}
\hskip25mm+u_0x\theta(z+b_3)+u_1xt_{1,3}\theta(z+b_1+b_3)+u_2xt_{2,3}\theta(z+b_2+b_3).
\end{equation}
This answer can be recovered by a straightforward (even though longer) computation using the gluings;
one can also easily see that it is correct by checking that it glues correctly and restricts to
the theta divisors on $\PP^1$ and $\PP^2$ that we computed previously. The computation of the vanishing
locus of the theta gradients is then also a straightforward combination of the previous results, and
we omit it.
One can also see in this case that the vanishing of the gradient of the theta function imposes $g$ independent conditions.

\subsection{Two $\PP^3$ bundles and an $F(2,2)$ bundle}
The only other type of semi-abelic varieties of torus rank 3 that form
a codimension $4$ locus in $\Perf$ (and also in $\Vor$) is the stratum corresponding in the toric picture to cutting up the cube
into an octahedron and two tetrahedra, see Figure \ref{fig:tetraocta} below.
\begin{figure}[ht!]
  \centering
    \includegraphics{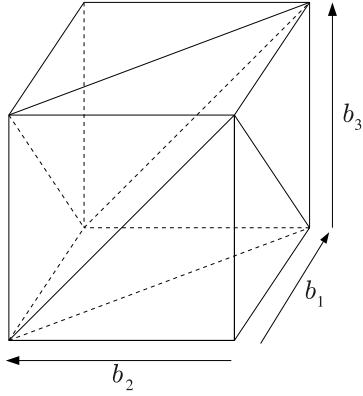}
\caption{Cube dissected into $2$ tetrahedra and $1$ octahedron}    \label{fig:tetraocta}
\end{figure}
The normalization of the toric variety has three irreducible components, two of which are $\PP^3$'s --- the toric varieties for which the corresponding polytope is a tetrahedron, while the toric variety corresponding to the octahedron is (see \cite{fultontoric})
$$
 F:=F(2,2)=\lbrace y_0y_1=y_2y_3=y_4y_5\rbrace\subset\PP^5,
$$
where we have denoted the homogeneous coordinates on $\PP^5$ by $y_0:y_1:y_2:y_3:y_4:y_5$.
This variety has $6$ nodes, they are the coordinate points of $\PP^5$.
In this description the eight faces of the octahedron correspond to the eight coordinate $\PP^2$'s contained in $F$,
corresponding to the case when all three products are equal to zero.
We will label such a $\PP^2$ by a triple $abc$ with
$a\in\lbrace 0,1\rbrace, b\in\lbrace 2,3\rbrace, c\in\lbrace 4,5\rbrace$, so that on this $\PP^2$ the
coordinates are $y_a:y_b:y_c$ (with the other three homogeneous coordinates being zero).

Here the gluing parameters
are crucial.
Indeed, the base abelian variety $B$ varies in ${\mathcal A}_{g-3}$, so we get $(g-3)(g-2)/2$ parameters.
We also have the three shifts $b_i \in B$ giving
another $3(g-3)$ parameters. Thus we must have two more free gluing parameters to make up codimension $4$.

\begin{figure}[ht!]
  \centering
    \includegraphics{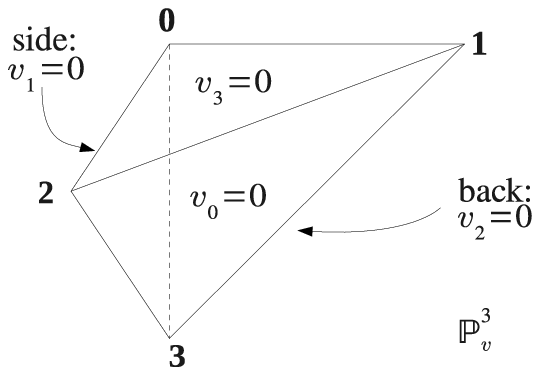}
    \includegraphics{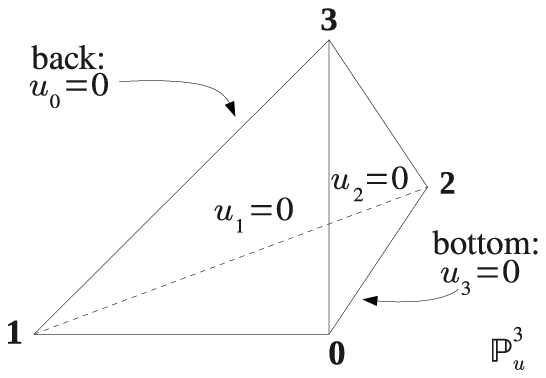}
\caption{Labeling of the tetrahedra} \label{fig:gluetetraocta_2}
\end{figure}

For the two $\PP^3$'s corresponding to the tetrahedra we choose homogenous coordinates $(u_0:u_1:u_2:u_3)$
and $(v_0:v_1:v_2:v_3)$ respectively. The faces of the tetrahedra are then given by
equations $u_i=0$ and $v_i=0$, respectively, for $i=0\ldots 3$.
Similarly to the case of two copies of $\PP^2$, we number the $u$ and the $v$ so as to label the
opposite sides with the same index, and for further use we also label the vertices by the number of the corresponding non-vanishing coordinate, see Figure \ref{fig:gluetetraocta_2}.

With a view towards further computations, similarly to the case of two copies of a $\PP^2$ bundle considered above, we
fix coordinates on $\PP^3_u$ and then choose the coordinates on $\PP^3_v$ to be their images under the
involution $j$, so that we have
\begin{equation}\label{jocta}
 j(z,(p:q:r:s))_u=(-b_1-b_2-b_3-z,(p:q:r:s))_v,
\end{equation}
which is easily seen to be compatible with the gluings by looking at the toric picture.

We now describe the semi-abelic theta divisors on these components.
Indeed, we know from the general theory what the restriction of the theta divisor to each fiber is (this
is given by the polytope description) and that the restriction to any section of the fibration
(such a section is isomorphic to a copy of $B$) is
a principal polarization on the base $B$ and hence the translate of a theta divisor.
In particular, on the fiber $F$ the theta divisor $T_F$ is a restriction of some
section of $\calO_{\PP^5}(1)$, while on each $\PP^3$ fiber the theta divisor is given by some
sections of $\calO_{\PP^3}(1)$.
It thus follows that there exist constants $\lambda_i,\mu_i,\nu_i\in\CC$ and points
$a_i,\alpha_i,\beta_i\in B$ such that the theta divisors are given by
$$
 T_F=\sum\limits_{i=0}^5 \lambda_i y_i\theta(z+a_i)|_F;\ \
 T_u=\sum\limits_{i=0}^3 \mu_i u_i\theta(z+\alpha_i); \ \
 T_v= \sum\limits_{i=0}^3 \nu_i v_i\theta(z+\beta_i).
$$
We now investigate the action of the involution, requiring it to fix the semi-abelic (reducible)
theta divisor. Writing down the condition that $j$ maps the zero locus of $T_u$ to
the zero locus of $T_v$ amounts to
$$
 \sum\limits_{i=0}^3\mu_iu_i\theta(z+\alpha_i)\sim
 \sum\limits_{i=0}^3\nu_iu_i\theta(-b_1-b_2-b_3-z+\beta_i)
$$
where $\sim$ denotes equality up to a non-zero constant factor. This implies
\begin{equation}\label{alphabeta}
\alpha_i=b_1+b_2+b_3-\beta_i
\end{equation}
for all the shifts, and the equality
$(\mu_0:\mu_1:\mu_2:\mu_3)=(\nu_0:\nu_1:\nu_2:\nu_3)$ as of projective points.
To simplify computations in what follows, we note that since the theta function is defined
only up to a constant multiple, we can choose $\mu_0=\nu_0=\lambda_0=1$.

The involution $j$ respects the fibration, i.e.~it lies over an involution
on the base $B$ and maps fibers to fibers. Each fiber over a fixed point of the involution
on $B$ is mapped to itself. We can see from the toric picture that on such a fiber the involution must act by $(y_0:y_1:y_2:y_3:y_4:y_5)\mapsto(t_1y_1:t_1^{-1}y_0:t_2y_3:t_2^{-1}y_2:t_3y_5:t_3^{-1}y_4)$
for some $t_1,t_2,t_3\in\CC^*$.
As in the case of torus rank $1$ degenerations,
we are again dealing with compactified $(\CC^*)^3$ bundles over the universal cover of the abelian variety $B$. Then the $t_i$ are sections of a bundle, namely the
trivial bundle, and thus constants.

Notice that we have the freedom of rescaling the
coordinates on $\PP^5\supset F$. We now choose coordinates in such a way as to get $t_1=t_2=t_3=1$ (note that
we could also have $-1$), so
that the involution on the normalization of the $F$ bundle over $B$ is then given by
\begin{equation}\label{jonf}
 j(z,(y_0:y_1:y_2:y_3:y_4:y_5))=
\end{equation}
$$
=(-b_1-b_2-b_3-z,(y_1:y_0:y_3:y_2:y_5:y_4))
$$

Writing down the condition $T_F\sim j^*T_F$
on the face $\PP^2\subset F$ given by $y_0=y_2=y_4$ yields
$$
 \lambda_1y_1\theta(z+a_1)+\lambda_3y_3\theta(z+a_3)+\lambda_5y_5\theta(z+a_5)\hskip30mm
$$
$$
 \sim y_1\theta(-b_1-b_2-b_3-z+a_0)+\lambda_2y_3\theta(-b_1-b_2-b_3-z+a_2)\hskip3mm
$$
$$
 \hskip49mm+\ \lambda_4y_5\theta(-b_1-b_2-b_3-z+a_4),
$$
which implies (after flipping the signs of the arguments of the even function $\theta$ on
the right-hand-side) in particular $a_1=b_1+b_2+b_3-a_0$ for the shifts, and for the coefficients we get
\begin{equation}\label{lambda}
(\lambda_1:\lambda_3:\lambda_5)=(1:\lambda_2:\lambda_4).
\end{equation}
By restricting $T_F$ to different coordinates $\PP^2$'s contained in
$F$ we get also $a_3=b_1+b_2+b_3-a_2$ and $a_5=b_1+b_2+b_3-a_4$, while all the equalities for the coefficients
together imply
$$
 (1:\lambda_1:\lambda_2:\lambda_3:\lambda_4:\lambda_5)= (\lambda_1:1:\lambda_3:\lambda_2:\lambda_5:\lambda_4)
$$
which, possibly after changing the signs of some coordinates
yields $\lambda_1=1, \lambda_3=\lambda_2, \lambda_5=\lambda_4$, so that finally we have
$$\begin{aligned}
 T_F&=y_0\theta(z+a_0)+y_1\theta(z+b_1+b_2+b_3-a_0)\\
 &+\lambda_2\left(y_2\theta(z+a_2)+y_3\theta(z+b_1+b_2+b_3-a_2)\right)\\
 &+\lambda_4\left(y_4\theta(z+a_4)+y_5\theta(z+b_1+b_2+b_3-a_4)\right).
 \end{aligned}
$$

\begin{figure}[ht!]
  \centering
    \includegraphics{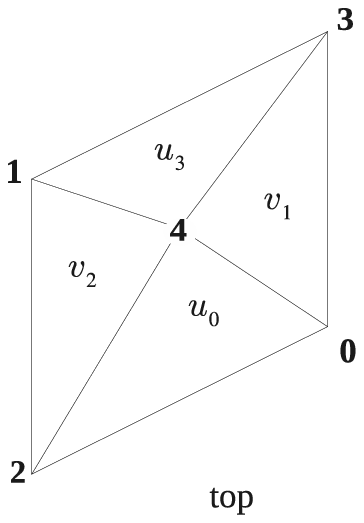}
    \includegraphics{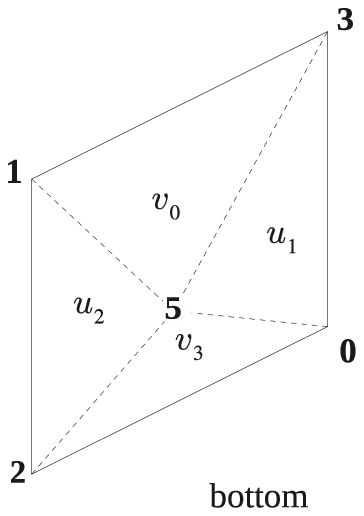}
\caption{The faces of the tetrahedron glued to the top and bottom of the octahedron}\label{fig:gluetetraocta}
\end{figure}

The shifts induced on the base $B$ of the semi-abelic variety under the gluings of the faces of
the cube are the same as for the gluing of the case $(\PP^1)^3$ discussed above: there
is no shift for gluings in the interior of the cube, and shifts by $b_i$ for gluing the faces,
which we label as in Figure \ref{fig:tetraocta}.
In Figure \ref{fig:gluetetraocta} we labeled the faces of the octahedron with the face of the
tetrahedron that gets attached to it ($u_i$ means that the coordinate $\PP^2$ given in $\PP^3_u$
by $\lbrace u_i=0\rbrace$ is attached).

For the gluing maps of the faces of the octahedron and of the tetrahedra, in addition to the shifts
we have to describe the permutation of the coordinates on the $\PP^2$'s being identified --- for this,
as before, we will see which coordinate points go where --- and the scaling parameters.
The vertices of the octahedron correspond to the 6 coordinate axes points $(0:\ldots:0:1:0:\ldots:0)\in\PP^5$
where only one homogeneous coordinate is non-zero.
In Figure \ref{fig:gluetetraocta} we label these points ${\bf 0},\ldots,{\bf 5}$; note that $({\bf 0},{\bf 1})$, $({\bf 2},{\bf 3})$,
and $({\bf 4},{\bf 5})$ are then the pairs of antipodal points on the octahedron.
To write down the gluing map from a face of the octahedron to the corresponding face of a
tetrahedron in coordinates, we need to make sure that the coordinate points get mapped appropriately.
For example, the ``$024$'' face of the octahedron (i.e.~the $\PP^2$
given by $\lbrace y_1=y_3=y_5=0\rbrace\subset F$) gets mapped to the $\lbrace u_0=0\rbrace$ face of the
first $\PP^3$, and looking at Figures \ref{fig:gluetetraocta_2} and \ref{fig:gluetetraocta}
we see that the axes points are
mapped ${\bf 0}\mapsto 1$, ${\bf 2}\mapsto 2$, and ${\bf 4}\mapsto 3$.
Thus the induced map on the $\PP^2$ must be given in coordinates by
$$
 (y_0:0:y_2:0:y_4:0)\mapsto (0:t_{00}y_0:t_{20}y_2:t_{40}y_4)_u,
$$
where $t_{00},t_{20},t_{40}\in\CC^*$ are the gluing parameters, and we use the subscript $u$ to keep
track of the fact that this is the first $\PP^3$. Since we are gluing projective planes,
the gluing is only defined up to scaling, and we can thus rescale to get
$t_{00}=1$. Furthermore, we can then also rescale coordinates $u_2$ and $u_3$ (without changing the coordinates $y_i$, which
we recall were fixed by prescribing the involution on $F$), to get $t_{20}=t_{40}=1$. Thus after these
choices we have fixed the coordinates $u_2,u_3$, and got the gluing to be
\begin{equation}\label{glue0}
 (z,(p:0:q:0:r:0))\mapsto (z,(0:p:q:r)_u)
\end{equation}
(recall that there is no shift on $B$ for gluing interior faces of a dicing of a cube); we have
switched to $(p:q:r)$ for coordinates on $\PP^2$ to make the other cases look similar.
The gluing of the $v_0=0$ face is obtained from this by acting by the  involution $j$, and is given by
$$
  (z,(0:p:0:q:0:r))\mapsto (z,(0:p:q:r)_v).
$$
Similarly, by combining the data of which faces get glued from Figure \ref{fig:gluetetraocta} with the depiction in Figure \ref{fig:tetraocta}
of the shifts on $B$ induced by the gluing, and the labeling of the vertices in Figure \ref{fig:gluetetraocta_2} (also shown in Figure 8, where we marked on the faces of the octahedron the labels of the vertices of the face of the tetrahedron glued to that face)
allows us to determine the other gluings. Looking at the fibers of the $\PP^2$ given by $u_1=0$, we have the map
$$
(y_0:0:0:y_3:0:y_5)\mapsto (t_{01}y_0:0:t_{21}y_5:t_{31}y_3)
$$
where we can choose $t_{21}$ to be equal to $1$, and then rescale the $u_0$ coordinate
so that $t_{01}=1$. Note that we have now fixed all the coordinates on $\PP^3_u$.
Including the base, the gluing on the
$\PP^2$ bundle given by $u_1=0$ is finally given by
\begin{equation}\label{glue1}
  (z,(p:0:0:q:0:r))\mapsto (z+b_1,(p:0:r:t_{31}q)_u).
\end{equation}
Similarly for the face $u_2=0$ we have
\begin{equation}\label{glue2}
  (z,(0:p:q:0:0:r))\mapsto (z+b_2,(t_{02}q:t_{12}r:0:p)_u),
\end{equation}
and for the gluing of the face $u_3=0$ we have
\begin{equation}\label{glue3}
  (z,(0:p:0:q:r:0))\mapsto (z+b_3,(r:t_{13}q:t_{23}p:0)_u).
\end{equation}
The gluings for the faces of the tetrahedron corresponding to $\PP^3_v$ are obtained from these by acting by the involution $j$.

We will now use the gluings to derive some identities between the parameters, ensuring that the restriction of $T_F$ to the $\PP^2$'s
contained in $F$  that are glued to the coordinate $\PP^2$s in the $\PP^3$'s agree with the restrictions of $T_u$ and $T_v$.
Indeed, for example we must have the restrictions
$$
 T_F|_{y_1=y_3=y_5=0}\qquad{\rm and}\qquad T_u|_{u_0=0}
$$
identified under the gluing map (\ref{glue0}). Writing out the formulae, this is equivalent to the relation
$$
 y_0\theta(z+a_0)+\lambda_2y_2\theta(z+a_2)+ \lambda_4y_4\theta(z+a_4)
\hskip4cm$$
$$ \hskip4cm \sim \mu_1y_0\theta(z+\alpha_1)+
 \mu_2y_2\theta(z+\alpha_2) +\mu_3y_4\theta(z+\alpha_3)
$$
being satisfied for all points $(y_0:y_2:y_4)\in\PP^2$ and for all $z\in B$. This in turn is equivalent to
having the identities
$$
 a_0=\alpha_1;\quad a_2=\alpha_2;\quad a_4=\alpha_3
$$
(as points of $B$) for the shifts, and the identity
$$
 (1:\lambda_2:\lambda_4)=(\mu_1:\mu_2:\mu_3)
$$
for the coefficients, as points in $\PP^2$.
Similarly from the other gluings we get further equalities relating $T_F$ and $T_u$. We first write down all the formulae
for the shifts on $B$, and will then deal with the gluing parameters. Indeed, from the gluings on the $u_1=0$ we get
$$
 a_0=\alpha_0+b_1;\ a_3=\alpha_3+b_1;\ a_5=\alpha_2+b_1,
$$
which using $a_3=b_1+b_2+b_3-a_2$ and $a_5=b_1+b_2+b_3-a_4$ and the above identities further yields
$$
 \alpha_1=\alpha_0+b_1;\qquad \alpha_2=b_2+b_3-\alpha_3.
$$
{}From the gluing on $u_2=0$ we then similarly get
$$
 b_1+b_3-\alpha_3=\alpha_1;\ \alpha_2=\alpha_0+b_2,
$$
so that expressing all the shifts in terms of $\alpha_0$ we obtain
\begin{equation}\label{shifts}
 \alpha_1=\alpha_0+b_1;\quad
  \alpha_2=\alpha_0+b_2;\quad
   \alpha_3=\alpha_0+b_3.
\end{equation}
Of course the gluings of the faces of $\PP^3_v$ will also give similar formulae
for the shifts: $\beta_i=\beta_0-b_i$. We can now write down two formulae for $a_2$:
in terms of $\alpha$, from the $u_0=0$ gluing, and in terms of $\beta$, from the $v_1=0$ gluing.
Using (\ref{alphabeta}) we then obtain
$$
 \alpha_0+b_2=\alpha_2=a_2=\beta_3-b_1=b_1+b_2+b_3-\alpha_3-b_1=
$$
$$
= b_2+b_3-(\alpha_0+b_3)=\alpha_0+b_2
$$
To summarize, we know at the moment that the semi-abelic theta divisor has on the
irreducible components the form
\begin{equation}\label{F}\begin{aligned}
 T_u&=u_0\theta(z)+\mu_1u_1\theta(z+b_1)+\mu_2u_2\theta(z+b_2)+\mu_3u_3\theta(z+b_3);\\
 T_v&=v_0\theta(z+b_1+b_2+b_3)+\mu_1v_1\theta(z+b_2+b_3)+\mu_2v_2\theta(z+b_1+b_3)\\
&\hskip79.5mm+\mu_3v_3\theta(z+b_1+b_2);\\
 T_F&=y_0\theta(z+b_1)+y_1\theta(z+b_2+b_3)+\lambda_2y_2\theta(z+b_2)\\
 &\hskip13.5mm+\lambda_2y_3\theta(z+b_1+b_3)+\lambda_4y_4\theta(z+b_3)+\lambda_4y_5\theta(z+b_1+b_2).
 \end{aligned}
\end{equation}

We now turn to determining the gluing parameters $t_{ij}$ in (\ref{glue1}), (\ref{glue2}), (\ref{glue3}) and the coefficients $\mu_i$
and $\lambda_i$ above.
Indeed, from (\ref{glue0}) we get from the
fact that $T_F|_{y_1=y_3=y_5=0}$ glues to $T_u|_{u_0=0}$ the identity
$$
 (1:\lambda_2:\lambda_4)=(\mu_1:\mu_2:\mu_3)
$$
(recall that $\lambda_0=\mu_0=1$).
{}From (\ref{glue1}) on $\lbrace y_1=y_2=y_4=0\rbrace\subset F$ glued to $\lbrace u_1=0\rbrace\subset \PP^3_s$ we get another identity
for $(1:\lambda_2:\lambda_4)$, this time involving a gluing parameter,
and altogether the identities we get from (\ref{glue0}),(\ref{glue1}),(\ref{glue2}),(\ref{glue3}) combined are
$$\begin{aligned}
(1:\lambda_2:\lambda_4)&=(\mu_1:\mu_2:\mu_3)=(1:t_{31}\mu_3:\mu_2)\\
&=(\mu_3:t_{02}:t_{12}\mu_1)=(t_{23}\mu_2:t_{13}\mu_1:1).
\end{aligned}
$$
Solving these, we get first $\mu_2=\lambda_4$, then
$$
 \mu_1=\lambda_2^{-1}\lambda_4;\quad \mu_3=\lambda_2^{-1}\lambda_4^2,
$$
and finally we compute
$$
 t_{31}=\lambda_2^{-2}\lambda_4^2;\quad t_{02}=\lambda_4^2;\quad t_{12}=\lambda_4^2;
 \quad t_{23}=\lambda_4^{-2};\quad t_{13}=\lambda_2^2\lambda_4^{-2}.
$$
Notice that as expected the two parameters $\lambda_2$ and $\lambda_4$ determine the
picture completely.

\smallskip
We now need to determine the fixed points of the involution, and compute the
gradients of the semi-abelic theta function there. Notice that in this case
for the first time we encounter a singular toric variety, with singularities
of $F$ being the coordinate points of $\PP^5$.
While considering the gradients of a function at a
singular point would be hard, we shall see that all fixed points lie either in the
interior part of $F$ (i.e. away from the planes where the gluing occurs) and thus
in the smooth locus, or on a $\PP^3$ corresponding to a tetrahedron, i.e. at least one of
the components of the normalization is smooth at these points. This will allow us
to use lemma \ref{lem:comparison}.

\smallskip
We first determine the fixed points of $j$ in the interior of $F$ (where none
of the $y_i$ are zero). As in all of the previous cases, note that our description
was of a trivial $F$ fibration over $\CC^{g-3}$, being a subset of $\PP^5\times\CC^{g-3}$.
Formula (\ref{jonf}) is for the involution there, and thus to determine the fixed
points on a bundle over $B$, we would need to include the exponential factors
similar to the ones needed above. These can be computed explicitly, but to
simplify the presentation (and since no new idea is necessary) we will omit
them in the following formulae. We will instead pretend that $\theta_m$ are even
or odd functions of $z$ depending on the characteristic (in fact they have
parity up to an exponential factor, which would precisely cancel the other
exponential factor from the transformation formulae for the $F$ bundle over $B$).
We thus determine the fixed points from
$$\begin{aligned}
 (z,(1:y_1:y_2:y_3:&y_4:y_5))=j (z,(1:y_1:y_2:y_3:y_4:y_5))\\
 &=(-b_1-b_2-b_3-z,(y_1:1:y_3:y_2:y_5:y_4)),
 \end{aligned}
$$
where we have rescaled the projective point with all non-zero coordinates to
get $y_0=1$. It follows that $z=m-\frac{b_1+b_2+b_3}{2}$ for some
two-torsion point $m\in B[2]$. For the projective coordinates we
must have $y_1=\pm 1$ and also $y_2=\pm y_3,$ $y_4=\pm y_5$. Since
we are only looking for fixed points on $F$, we finally see that
$y_1,y_3,y_5\in\lbrace -1,1\rbrace$ and finally that the
fixed points are of the form
$$
 \left(m-\frac{b_1+b_2+b_3}{2},(1:y_1:y_1y_3:y_3:y_1y_5:y_5)\right)
$$
for $y_1,y_3,y_5\in\lbrace -1,1\rbrace$. The number of such fixed points is thus
equal to $2^3\cdot 2^{2(g-3)}$.

Evaluating the semi-abelic theta divisor at such a point gives (up to
the exponential factors which we omitted in the expressions for $y_i$, and
which cancel with the exponential factors resulting from switching
to $\theta_m$)
$$
 \theta_m\left(\frac{b_1-b_2-b_3}{2}\right)+y_1\theta_m\left(\frac{b_2+b_3-b_1}{2}
 \right)
$$
$$+\lambda_2y_3\left(y_1\theta_m\left(\frac{b_2-b_1-b_3}{2}\right)+
  \theta_m\left(\frac{b_1+b_3-b_2}{2}\right)\right)
$$
$$
+\lambda_4y_5\left(y_1\theta_m\left(\frac{b_3-b_1-b_2}{2}\right)+
  \theta_m\left(\frac{b_1+b_2-b_3}{2}\right)\right),
$$
which, as expected, vanishes for $m$ odd and $y_1=1$ or for $m$ even and $y_1=-1$.

It remains to compute the gradient of the theta function at such a fixed point. These two cases are similar, and we give the formulae for the case of $m$ odd and $y_1=1$.
Note that we are working on $F$
and that we can use $z$, $y_1$, $y_3$, $y_5$ as local coordinates. For computing the gradient of $T_F$ near a fixed point as above we locally set $y_0=1$, and write $y_2=y_1/y_3$, $y_4=y_1/y_5$, so that we get for the gradients evaluated at these points
$$
 \frac{\partial T_F}{\partial y_1}= \theta_m\!\left(\frac{b_2+b_3-b_1}{2}\right)\!
 -\lambda_2y_3\theta_m\!\left(\frac{b_1+b_3-b_2}{2}\right)\!
 -\lambda_4y_5\theta_m\!\left(\frac{b_1+b_2-b_3}{2}\right)\! ;
$$
$$
  \frac{\partial T_F}{\partial y_3}= 2\lambda_2\theta_m\left(\frac{b_1+b_3-b_2}{2}\right);
$$
$$
  \frac{\partial T_F}{\partial y_5}= 2\lambda_4\theta_m\left(\frac{b_1+b_2-b_3}{2}\right).
$$
(so that these three all vanish if and only if the three values $\theta_m((b_1+b_2+b_3)/2-b_i)$ all vanish), and
$$
 \frac{\partial T_F}{\partial z}= 2\frac{\partial \theta_m}{\partial z}\left(\frac{b_1-b_2-b_3}{2}\right)+2y_3\lambda_2
 \frac{\partial \theta_m}{\partial z}\left(\frac{b_2-b_1-b_3}{2}\right)
 +
$$
$$
 + 2y_5\lambda_4\frac{\partial \theta_m}{\partial z}\left(\frac{b_3-b_1-b_2}{2}\right).
$$
Studying the dimension of the locus where these gradients are all zero simultaneously
is complicated (the $z$-derivative is hard to deal with).
However, the stratum of semi-abelic varieties of this type is already codimension 4 in $\Perf$, and it is enough for our purposes to show that the common zero locus of these partial derivatives
has codimension at least $2$ on this stratum.
Since the vanishing of the $y_1,y_3,y_5$ derivatives above implies that three independent points lie on the theta divisor, these vanishing conditions are in codimension at least $3$ within the locus of semi-abelic varieties with such a toric structure.

\smallskip
We now need to deal with the fixed points of the involution $j$ that are on the ``gluing'' locus,
i.e.~where the map from the normalization to the semi-abelic variety is not $1$-to-$1$.
A geometric picture of the dicing of the cube as shown in Figure \ref{fig:tetraocta} (or a direct calculation)
shows that there are no such
fixed points on the faces of the tetrahedra or of the octahedron, and that the only fixed points are on the
edges or at the vertices. Note now that all edges of the octahedron are glued to some edges of the
tetrahedra (and the total number of edges after the gluings is $6$), while all the vertices are glued together.
We will now compute the fixed points on the tetrahedra, and for gradients
will work with $T_u$ and $T_v$ rather than $T_F$.
For example the fixed points on the ``edge'' $u_0=u_1=0\ne u_2u_3$ are computed from the
equation
$$
 (z,(0:0:1:p)_u)=j(z,(0:0:1:p)_u)
$$
$$
=(-b_1-b_2-b_3-z,(0:0:1:p)_v) \sim
 (-b_1-b_2-b_3-z,(0:0:0:1:0:p))
$$
$$
\sim(-b_2-b_3-z,(0:0:p:t_{31})_u),
$$
which gives the fixed points as
$$
 \left(m-\frac{b_2+b_3}{2},(0:0:1:\pm t_{31}^{1/2})_u\right).
$$

Evaluating the semi-abelic theta divisor at these points gives (up to the exponential factors)
$$
 T_u=\mu_2\theta_m\left(\frac{b_2-b_3}{2}\right)\pm \mu_3 t_{31}^{1/2}
 \theta_m\left(\frac{b_2-b_3}{2}\right),
$$
which upon using $\mu_2=\lambda_4, \mu_3=\lambda_2^{-1}\lambda_4^2,
t_{31}=\lambda_2^{-2}\lambda_4^2$ becomes simply
$$
 \lambda_4\theta_m\left(\frac{b_2-b_3}{2}\right)\pm\lambda_4\theta_m\left(\frac{b_3-b_2}{2}\right).
$$
Thus such a point generically lies on the semi-abelic theta divisor if
$m$ is odd and the $+$ sign is chosen or if $m$ is even and the $-$ sign is chosen.
For the gradient of $T_u$ at such a point we then get
$$
 \frac{\partial T_u}{\partial u_0}= \theta_m\left(\frac{-b_2-b_3}{2}\right);
 \quad
  \frac{\partial T_u}{\partial u_1}=\mu_1\theta_m\left(\frac{b_1-b_2-b_3}{2}\right);
$$
$$
  \frac{\partial T_u}{\partial u_2}=\mu_2\theta_m\left(\frac{b_2-b_3}{2}\right),
$$
which clearly imposes at least $3$ independent conditions.

The computation for the other edges is completely analogous. The total number of fixed
points of $j$ of this type is then equal to $6\cdot 2\cdot 2^{2(g-3)}$, each of them
being the limit of $4$ different $2$-torsion points on smooth ppavs.

\smallskip
Finally we need to compute the fixed points of $j$ over the ``vertices'' of the dicing.
These are computed from the equation
$$\begin{aligned}
 (z,(1:0:0:&0)_u)=j(z,(1:0:0:0)_u)\\
&=(-b_1-b_2-b_3-z,(1:0:0:0)_v)\\
&\sim(-b_1-b_2-z,(0:0:0:0:0:1))\\
& \sim(-b_1-b_2-z,(0:0:0:1)_v)\\
&\sim(-b_1-z,(1:0:0:0:0:0))
 \sim(-z,(1:0:0:0)_u),
 \end{aligned}
$$
which yields simply the fact that $m$ is a 2-torsion point on $B$. Since all the vertices are glued together, we thus have $2^{2(g-3)}$ fixed points, each arising
as the limit of 8 two-torsion points on smooth ppavs. For a consistency check, we have the valid equality
$$
 2^{2g}=2^{2(g-3)}(8+6\cdot2\cdot 4+8)
$$
for the number of limits of two-torsion points on smooth ppavs. The gradient computation is completely analogous to the previous cases and we omit it.

\section{The codimension 5 strata of non-standard semi-abelic varieties in $\Perf$}
In this section we explain the approaches to handling all the non-standard degeneration types
that give strata of codimension 5 in $\Perf$. Instead of describing, very laboriously as above,
the semi-abelic theta divisors explicitly, we note that our main theorem \ref{theo:Gg} states essentially that $\overline{G^{(g)}}$ does
not have codimension 5 components contained in the boundary. Thus it
suffices for our purposes to show that none of these codimension 5 strata can be contained in
$\overline{G^{(g)}}$. To do this it suffices to exhibit one point not contained in the
corresponding variety, and this is done by degeneration. Below we give two techniques to do this, one that
works easily when the toric variety is a product, and thus handles many (but not all) of our cases, and the
second, by degenerating to the ``principal'' degenerations --- dicing into simplices --- that applies in all
of our situations. The two approaches seem to be independently interesting and general.

\subsection{Approach via ``semi-direct'' products}
To handle the first two non-standard codimension 5 cases in $\beta_4^0$, where the fibers of the normalizations of the semi-abelic varieties are
of the form $\PP^1 \times (F(2,2) \sqcup 2 \PP^3)$ or $\PP^1\times 2(\PP^1\times\PP^2)$ we will use the corresponding cases of $F(2,2)
\sqcup 2\PP^3$ and $2(\PP^1\times\PP^2)$ handled above.
We shall explain this for the first of these two cases, the second being analogous --- the main point is to use the $\PP^1$ product
structure.
We start with a semi-abelic variety of torus rank $1$, i.e., a semi-abelic variety whose normalization
is a $\PP^1$ bundle over an abelian variety $A$ of dimension $g-1$. Now we choose a degeneration of $A$ to a semi-abelic
variety $X$ of torus rank $3$ with base an abelian variety $B$ of dimension $g-4$,
where the fiber of the normalization of $X$ is of the form
$(F(2,2) \sqcup 2 \PP^3)$.
Under such a degeneration, the divisor of the polarization behaves
as follows. For a semi-abelic variety of torus rank 1 we start with a function of the form $\theta(z) + x \theta(z+b)$ where $z \in A$, $x$ is the coordinate
of $\PP^1$ and $b \in A$ is the shift parameter.
Degenerating this in the family described above we obtain on each component of the normalization a function of the
form $T(z') + xt_{b'}^*(T(z'))$ where $z'$ are coordinates on the total space of an $F(2,2)$ or $\PP^3$ bundle over $B$ and $b'$ is
a point of the semi-abelian variety which is a rank $3$ extension of $B'$ and which acts by $t_{b'}$ on the semi-abelic variety (and thus
also on the components of the normalization).  We have already seen that it is a non-empty condition that the theta divisor on $X$ has a singular point
at a fixed point of the involution.
We can then argue as in the discussion of formulae (\ref{equ:partial1}) or
(\ref{equ:partial2}) in section \ref{sec:thetagradients} that lying in
$\overline{G_g}$ imposes a non-trivial condition for the $\PP^1$ bundles
over $X$ considered here. This is sufficient for us since the loci under consideration (being proper subvarieties of $\beta_4$) already have codimension $5$ in $\Perf$.

\subsection{Approach via degenerating to a principal semi-abelic variety}
An approach to handle all the non-standard codimension 5 strata is by dicing the cube further, into simplices. Note that all our non-standard codimension 5 strata correspond to cones spanned by forms $x_i^2$ or $(x_i-x_j)^2$. These cones are thus contained in the corresponding torus rank $k$ {\it principal
cone} $\sigma_0=\langle x_i^2, (x_i - x_j)^2 | \, i,j= 1, \ldots k \rangle$, the
so called {\it {principal cone}}. Thus the strata of semi-abelic varieties corresponding to each of these cases
contain the locus of semi-abelic varieties corresponding to the principal cone. The principal cone defines in $\calA_g^{\operatorname{Perf}}$ a stratum with the normalization of the corresponding semi-abelic variety being $k!$ copies of a $\PP^k$ bundle over some $B\in\calA_{g-k}$. We call these semi-abelic varieties the {\it {principal}} semi-abelic varieties
of torus rank $k$ (For $k=4$ there is only one other dimension $10$ cone in the perfect cone decomposition,
namely the so called second
perfect cone, which, however, is not a cone in the second Voronoi decomposition where it gets subdivided.)
\begin{prop}\label{prop:principalok}
The stratum of principal semi-abelic varieties of torus rank $k$ is not contained in
$\overline{G^{(g)}}$.
\end{prop}
\begin{proof}
To prove this proposition, we will first describe explicitly the semi-abelic theta divisor on the principal
semi-abelic varieties, then will discuss the gradients at the fixed points of $j$ to deal with $\overline{G^{(g)}}$.

There are finitely many combinatorial choices for how a $k$-dimensional cube is diced into $k$-dimensional simplices; the
resulting semi-abelic theta divisor would depend on this, but will be very similar for all the cases.
Since we have the freedom of choosing coordinates on half of the $\PP^k$'s
independently (the other half would be given by the involution), by rescaling the coordinates we can ensure that on
each $\PP^k$ the semi-abelic theta divisor $T_i$ (for $i=1\dots k!$) has the form
$$
 T_i(z,(u_0:\dots:u_k))=\sum\limits_{n=0}^k u_n\theta(z+v_n^i)
$$
for some shifts $v_n^i\in B$. The gluings all happen along the faces of the simplices, i.e.~along $\PP^{k-1}$'s, and there
are no free scaling parameters for the gluings, while the shifts are either $0$ for the gluings in the interior or
the appropriate $b_i$ for gluing the exterior faces of the cube (each subdivided into $(k-1)!$ simplices). The gluing
conditions for $T_i$ and $T_j$ on some face would then be to say that all, possibly except one (for the variable the vanishing
of which gives the faces of the simplex that are on the faces of the cube) shifts $v_n^i$ are obtained from the shifts $v_n^j$ by
a suitable permutation (which tells us how the coordinates on the glued $\PP^{k-1}$ faces of $\PP^k_i$ and $\PP^k_j$
are identified), while if we are gluing faces of the simplices that are faces of the cube, we shift by some $b_i$.

We start from the simplex at the origin, for which $k$ of its
faces are glued to other simplices by $+b_i$ shifts, and the one face in the interior of the cube is glued without the shift. The
``origin'' of the cube is a vertex of this simplex that is identified with each of the other $k$ vertices of the simplex
via a shift by the corresponding $b_i$. It thus follows that $\theta(z+v_0^1)$ (the value of $T_1$ for $u_1=\dots=u_k=0$)
is identified with $\theta(z+v_n^1)$ via a shift by $b_n$, and it thus follows that $v_n^1=v_0^1+b_n$. By an overall
coordinate shift on $\PP^k$ we can then make
$$
  T_1(z,(u_0:\dots:u_k))=u_0\theta(z)+\sum\limits_{n=1}^k u_n\theta(z+b_n).
$$
We note that in general for any simplex all of its vertices are vertices of the cube, and they are
all identified together by certain shifts. Since there is a non-zero shift for identifying any vertex of
the cube to any other vertex, none of these shifts are zero, and thus it follows that $v_n^i\ne v_m^i$ for
$n\ne m$.

This consideration determines the semi-abelic theta divisor on any $\PP^k_i$ up to an overall shift,
and this shift can then be obtained by writing down the gluing conditions to eventually glue to $\PP^k_1$. The answer
for $k=2$ is discussed in detail in the section on the non-standard torus rank 2 degeneration. As an example of
what happens for higher torus tank, here is a possible answer for $k=3$ (depending on how the cube is
diced, the other symmetric choices may occur):
$$
 v^1=(0,b_1,b_2,b_3);\quad v^2=(b_1+b_2,b_1,b_2,b_3);\quad v^3=(b_1+b_2,b_1,b_1+b_3,b_3);
$$
$$
 v^6=(b_1+b_2+b_3,b_2+b_3,b_1+b_3,b_1+b_2);\quad v^5=(b_3,b_2+b_3,b_1+b_3,b_1+b_2);
$$
$$
 v^4=(b_3,b_2+b_3,b_2,b_1+b_2).
$$

Now in order to investigate the fixed points of the involution $j$, we first note that since it permutes the simplices,
all the fixed points are on the gluing loci. To determine the gradient of the semi-abelic theta function at each such
point, one can take any of the components of the normalization containing it. The partial derivative of $T_j$ with respect to $u_n$ is equal
to $\theta(z+v_n^j)$. Note, however, that for $b_i$ generic no two points $v_n^j$ differ by a point of order two, and thus
we cannot have both $z+v_n^j$ and some $z+v_m^j$ be points of order two on $B$. Thus generically some $\theta(z+v_n^j)$
does not vanish, and thus a generic principal semi-abelic variety does not lie in $\overline{G^{(g)}}$.
\end{proof}
\begin{rem}
One could try to prove the same statement about $N_2$: that it does not contain the stratum of principal degenerations. To do this, one would
use the setup and the results of Ciliberto and van der Geer \cite{civdg2}.
By \cite[Rem.~18.1]{civdg2}, the scheme of vertical singularities is defined over $\Perf-\beta_5-Z$, (where $Z$ is the locus where the perfect cone and the Voronoi compactifications differ)
and can in fact be defined for any principal degeneration (note that this stratum belongs to both the perfect cone
and the Voronoi compactification).
Similarly to the discussion of the non-standard (i.e.~principal) torus rank 2 degenerations in
\cite[Sec.~17]{civdg2}, the vertical singularities are of different types, depending on how many of the homogeneous
coordinates on the corresponding $\PP^k_j$ vanish. Indeed, the condition for there to exist $(u_0:\dots:u_k)$ with
$u_0\dots u_k\ne 0$ such that the $T_j$ is singular at that point is for the $k+1$ theta divisors $\Theta_{v_n^j}$
to be tangentially degenerate, in the notation of \cite[Sec.~11]{civdg2}, which is to say that for some $z\in B$
the Gauss images of the $k+1$ points $\theta(z+v_n^j)$ are linearly dependent (or for all the translates of the theta divisor to be singular, or a limit of these two cases).

To complete the argument, one would need to apply a suitable version of
\cite[Prop.~12.1]{civdg2}. Unfortunately, the proof given there seems incomplete. The $k=1$ case of the proposition, the only one used in that paper, holds, but it seems that the argument given only proves a weaker statement for higher $k$, insufficient to complete our proof. We have not been able to adapt the argument to our situation.
\end{rem}
As a corollary we immediately obtain
\begin{cor}\label{cor:codimension5ok}
The locus
$\overline{G^{(g)}}$ contains no strata of non-standard codimension $5$ consisting of non-standard semi-abelic varieties.
\end{cor}
\begin{proof}
Since the corresponding cones of these strata are all contained in the principal cone it follows that the closure of the
codimension $5$ strata contain the locus of principal semi-abelic varieties of torus rank $4$. Thus the assertion follows
immediately from proposition \ref{prop:principalok}.
\end{proof}

\subsection{Explicit formulae for the case of two $\PP^3$ bundles and two quadric cone bundles}
The two approaches above allow one to deal with all the types of non-standard semi-abelic varieties
giving strata of codimension 5 in $\Perf$ without having to undertake a laborious computation similarly
to the one for the $F(2,2)\sqcup 2\PP^3$ toric part. In fact, if the toric part of a semi-abelic variety is
a product of $\PP^1$ and another known case, the formulae for the semi-abelic theta divisor can be obtained
rather straightforwardly.

For completeness, in this section we give the explicit expressions for the involution, gluings, and the semi-abelic
theta divisor for the torus rank 3 semi-abelic variety corresponding to the case of two tetrahedra and two square
pyramids. This will therefore complete the explicit computations for {\it all types} of semi-abelic varieties of torus rank up to 3.
The derivations in this case are tedious, and for brevity we will omit all of them, just explaining
the coordinate choices made, and writing down the resulting formulae.

\smallskip
As in the previous cases, we will choose coordinates on $\PP^3_v$ to be the image
of those on $\PP^3_u$ under the involution $j$, and similarly the coordinates on $\PP^4_x$ (where one
quadric cone sits)
to be the image of those on $\PP^4_y$ (where the other quadric cone sits)
under the involution, while we choose the coordinates on the bases
of these bundles in such a way that as for the octahedron case we have $j:z\mapsto
-b_1-b_2-b_3-z$, as in (\ref{jocta}). Writing down the general expressions for theta divisors on each $\PP^3$
(where it is a section of $\calO(1)$ on the fiber), and on each quadric cone (where it is a restriction
of a section of $\calO(1)$ from $\PP^4$), then suitably changing coordinates on $\PP^3_u$ and on $\PP^4_x$,
while preserving the cone, and then requiring the gluings and the involution to preserve the theta divisor allows to
compute all the unknown coefficients in terms of one free parameter that we denote $c\in\CC^*$.

The gluings are then given by
$$\begin{aligned}
 (z,(p:q:0:r:0))_x&\mapsto (z+b_3,(p:q:r:0))_u;\\
 (z,(p:q:0:0:r))_x&\mapsto (z-b_2,(q:p:0:r))_v;\\
 (z,(p:0:q:r:0))_x&\mapsto (z-b_1,(c^2r:0:p:q))_v;\\
 (z,(p:0:q:0:r))_x&\mapsto (z,(0:r:c^{-2}q:p))_u,
 \end{aligned}
$$
for the identifications of the faces of the tetrahedra with the sides of the pyramids (along $\PP^2$'s). The new feature of this case compared to that of a dicing into an octahedron and two tetrahedra is the
gluings of the bases $x_0=0$ and $y_0=0$ of the two pyramids (along $\PP^1\times\PP^1$), which is given by
\begin{equation}\label{gluebase}
 (z,(0:p:q:r:s)_x)\mapsto (z,(0:c^{-1}q:cp:c^{-1}s:cr)_y).
\end{equation}
The resulting expressions for the components of the semi-abelic theta divisors are
$$\begin{aligned}
 T_u&=u_0\theta(z)+u_1\theta(z+b_1)+cu_2\theta(z+b_2)+u_3\theta(z+b_3);\\
 T_x&=x_0\theta(z+b_3)+x_1\theta(z+b_1+b_3) +c^{-1}x_2\theta(z+b_2)\\
&\hskip24mm+ cx_3\theta(z+b_2+b_3)
 +x_4\theta(z+b_1),
 \end{aligned}
$$
and $T_v$ and $T_y$ are obtained from these by applying the involution $j$.

\section{Proofs of the main theorems}
In this section we use the explicit descriptions of the semi-abelic theta divisors given above to obtain the proofs of our main theorems.

\begin{rem}
To attempt to prove conjecture \ref{theo:N2}, one would use the method
pioneered in \cite{civdg2}: by studying the intersection of the closure $\overline{N_2}$ with the boundary strata. Note that since we have described more boundary strata, we have more control of the situation, and are thus able to push the
computation deeper into the boundary. Recall from \cite[Rem.~18.1]{civdg2} that away from $\beta_5$ (and also for the standard torus rank
5 degenerations) the vertical singularities of the semi-abelic theta divisors can be defined geometrically. We recall that the way this is done is by using the log-differentials $\Omega^1_\Perf(\log D)$, which is to say that one computes the gradient of the semi-abelic theta divisor with respect to the coordinate $z\in B$ on the base of the semi-abelic variety, and to the coordinates on the fiber, and then has to take into account the ``gluing'' singularities, see eg.~\cite[formula (28)]{civdg2} for the case when the toric variety is $2\PP^2$.

One would now proceed
by studying $\overline{N_2}\cap\partial\Perf$.
Indeed, assume for contradiction that $N_2$ had an irreducible component of codimension 4 in $\calA_g$ (the results of Ciliberto and van der Geer \cite{civdg2} apply to show there could not be a component of codimension 3). By the result of Keel and Sadun \cite{kesa} this component would have to intersect the boundary of $\Perf$, and we let $M$ be its intersection with the boundary, which is thus of codimension 5 in $\Perf$.
If $M$ were to intersect the locus of standard semi-abelic varieties, the vertical singularities there are completely described in \cite{civdg2}, and noting that an analog of lemma 16.1 there applies to show that if the dimension of the space of vertical singularities is at least 2, $B\in N_1$,
we see that $M$ cannot have a codimension 5 component intersecting the locus of standard semi-abelic varieties. The situation with non-standard torus rank 2 semi-abelic varieties is also handled in \cite{civdg2}, and thus it remains for us to deal with the possibility that $M\subset\beta_3$; note also that $M$ cannot contain $\beta_5$ (where the set of standard compactifications of torus rank 5 is open), and thus $M$ must intersect $\beta_3^0\sqcup\beta_4^0$, and moreover it must intersect some locus of non-standard compactifications there. For all codimension 5 loci corollary \ref{cor:codimension5ok} applies, and thus it remains to handle the two loci of non-standard torus rank 3 compactifications that have codimension 4 in $\Perf$ --- when the toric part is $2\PP^1\times\PP^2$ or $F_{2,2}\sqcup 2\PP^3$.

The case of two copies of $\PP^1\times\PP^2$ is straightforward --- one can use the description of the semi-abelic theta divisor
given by (\ref{p1p2}) and compute directly, or note that this case arises as the limit of the standard torus rank 3 case when
one of the gluing parameters approaches zero (this of course also equivalent to the case when it approaches infinity),
so that all the equations and descriptions arise as limits of those for the standard case.

The hardest case to deal with is perhaps that of the toric part being $F(2,2)\sqcup 2\PP^3$. Note that here the normalization of the toric variety is itself singular.
Here one would need to have an appropriately corrected and generalized version of \cite[Prop.~12.1]{civdg2}, and also a geometric understanding of the singularities of the theta divisors possibly appearing at singular points of the $F(2,2)$ fibration. This seems plausible, but, as before, we were not able to complete this step.
\end{rem}

\begin{proof}[Proof of theorem \ref{theo:Ig}]
We want to (re)prove the statement that the codimension of $I^{(g)}$ in $\calA_g$ is equal to $g$ for $g\le 5$. Indeed, for $g=1$ this is trivial. For $g=2$ the theta divisor of an indecomposable ppav is smooth, and the theta divisor of a product of two elliptic has one singular point which is not of multiplicity three. For higher $g$, suppose the codimension of some irreducible component $Z$ of $I^{(g)}$ is equal to $d$. We then know, since $I^{(g)}$ is given locally by $g$ equations, that $d\le g$, and we want to prove the equality. For $g\ge 3$, by a theorem of Keel and Sadun \cite{kesa} any codimension $g$ subvariety of $\calA_g$ is non-complete, and thus the closure $\overline Z$ of $Z$ in $\Perf$ must intersect the boundary of $\Perf$. Moreover, the codimension of $Y:=\overline Z\cap\partial\Perf$ within $\Perf$ must then be equal to $d+1$. If $g\le 5$ (which are the only cases when we can prove the theorem), if $d<g$, it would follow that the codimension of $Y$ in $\Perf$ is at most 5, and thus $Y$ must intersect one of the strata of semi-abelic varieties of types described above.
For each of these cases from the description of the vanishing locus of $\tilde f_m$ (given at the end of the corresponding section of the text) we see that if theorem \ref{theo:Ig} holds in genus $g-k$, where $k$ is the torus rank, then the intersection of $\overline{G^{(g)}}$ (and thus of $Y$ contained in it) with the corresponding stratum is codimension strictly greater than $g$ within $\Perf$. It thus follows that the codimension of $\overline Z$ in $\Perf$ is at least $g$.
\end{proof}

\begin{proof}[Proof of theorem \ref{theo:Gg}]
Recall that this theorem says that the vanishing locus of a function $\tilde f_m$ on $\Perf$, which is the
extension to the boundary of the gradient of the theta function at the corresponding odd two-torsion point,
does not have an irreducible component of codimension less than 6 contained in the boundary.
Indeed, in the previous sections we described explicitly the functions $\tilde f_m$ for each type of
semi-abelic varieties giving a stratum of codimension at most 5 in $\Perf$. For the simplest case of
torus rank 1 degenerations we proved proposition \ref{codim1Ok} showing that $\overline{G^{(g)}}$, the zero
locus of $\tilde f_m$, does not have any irreducible components contained in $\beta_1^0$. In the following sections
above we then handled all the types of semi-abelic varieties in $\beta_2^0\sqcup\beta_3^0\sqcup\beta_4^0$ that
give loci of codimension at most 5 in $\Perf$. From the explicit descriptions of the theta divisors, we know,
as in the proof of proposition \ref{codim1Ok}, that the intersection of $\overline{G^{(g)}}$ with the
corresponding locus is of codimension greater
than $g$ in $\Perf$. On the other hand, since the zero locus of the gradient $\tilde f_m$ is locally
given by the vanishing of $g$ functions, any irreducible component of $\overline{G^{(g)}}$ has codimension at
most $g$ in $\Perf$, and we have arrived at a contradiction. It thus follows that $\overline{G^{(g)}}$ has
no irreducible components of codimension less than 6 in $\Perf$ that are contained in the boundary of $\Perf$,
and the theorem is thus proven.
\end{proof}

\end{document}